\documentclass{article}

\usepackage[a4paper,top=2cm,bottom=2cm,left=3cm,right=3cm,marginparwidth=1.75cm]{geometry}

\usepackage{graphicx}
\usepackage{amsmath}

\DeclareMathOperator*{\argmin}{arg\,min}

\usepackage[colorlinks=true, allcolors=blue]{hyperref}
\usepackage[ruled,vlined]{algorithm2e}
\usepackage{amsfonts}
\usepackage{multicol}

\usepackage[normalem]{ulem}

\newcommand{\id}{\,\text{d}}
\newcommand{\dpfrac}[2]{\dfrac{\partial #1}{\partial #2}}
\newcommand{\dint}{\displaystyle\int}
\newcommand{\dsum}{\displaystyle\sum}
\renewcommand{\vec}[1]{\mathbf{#1}}
\newcommand{\tf}[1]{\bar{#1}}
\newcommand{\ts}[1]{\bar{\mathcal{#1}}}
\newcommand{\vx}[0]{\mathbf{x}}
\newcommand{\onep}[0]{\emph{one-phase}}
\newcommand{\twop}[0]{\emph{two-phase}}

\newcommand\xm{%
  X\kern-.0em%
  --%
   \hbox{M}\kern+.16em%
   \raise-0.65ex\hbox{--}\kern-.7em%
  \raise0.34ex\hbox{es}%
  \hbox{{h}}%
\@}

\def\iGamma{\rule{0pt}{1.5ex}\Gamma}

\usepackage{authblk}

\title{Solving the Porous Medium Equation with the eXtreme Mesh deformation approach (\xm)}
\author[1]{Alexandre {Chemin}}
\author[1]{Jonathan {Lambrechts}}
\author[1]{Nicolas {Moës}}
\author[1]{Jean-François {Remacle}}
\affil[1]{Institute of Mechanics, Materials and Civil Engineering (iMMC), Avenue Georges Lemaître 4, 1348 Louvain-la-Neuve, Belgium}
\date{November 2024}

\providecommand{\keywords}[1]{\textbf{\textit{Keywords: }} #1}

\begin{document}
\maketitle

\begin{abstract}
  We introduce a new scheme for solving the non-regularized Porous Medium Equation. It is mass conserving and uses only positive unknown values. To address these typically conflicting features, we employ the eXtreme Mesh deformation approach (\xm), specifically designed for problems involving sharp interfaces. The method ensures that the interface is always meshed, even in the face of complex topological changes, without the need for remeshing or altering the mesh topology. We illustrate the effectiveness of the approach through various numerical experiments.
\end{abstract}


\section{Introduction}

The Porous Medium Equation (PME) is a non-linear partial differential equation of significant interest due to its widespread applicability in modeling various physical phenomena such as gas ﬂow in porous medium, incompressible ﬂuid dynamics and non-linear heat transfer \cite{hayes2012physics, leibenzon1947motion, aronson2006porous, vazquez2007porous}. 

The PME is usually written as the following initial value problem:
\begin{eqnarray}
  \frac{\partial u}{\partial t} &=& \nabla \cdot \left(\kappa u^m \nabla u\right),
    \label{eq:PME}  \\
  u({\vec x},t_0) &=& u_0({\vec x})  \geq 0.
    \label{eq:PMEINIT}
\end{eqnarray}

Here, $m > 0$ is a real parameter, $u({\vec x},t)$ is the unknown scalar function in space $\vec x \in  \Omega \subset\mathbb{R}^d$ and time $t \in [t_0, T]$ and $\kappa u^m$ is homogeneous to a diffusivity coefficient. The case $m<0$ corresponds to the Fast Diffusion Equation and is not treated in this paper. In Eq.~\eqref{eq:PME},  $u$ is frequently employed to denote physical quantities such as gas density or temperature, requiring the preservation of non-negativity in solutions.
The case $m = 0$ leads to the standard heat equation.

The PME is parabolic provided $u \neq 0$, and degenerates when $u = 0$ \cite{vazquez2007porous}. The nature of the solutions to Eq.~\eqref{eq:PME} differs fundamentally depending on whether we consider $m=0$ (linear case) or $m>0$ (non-linear case). In the non-linear scenario, if the initial distribution $u_0({\vec x})$ of the quantity of interest $u$ has a compact support, then the solution $u(x,t)$ will have a compact support for $t>0$. The region $\mathcal{Q}(t)\subset\Omega$ where the solution is zero will shrink in size over time but remain present. Thus, in the non-linear case, a sharp interface $\Gamma(u)$ may exist, separating $\mathcal{Q}(t)$ from the rest of the domain where the solution is positive.  This behavior arises from the finite propagation PME's property which is defined in \cite{vazquez2007porous} as: \emph{disturbances from the level $u=0$ propagate in time with finite speed for solutions of the PME}. This constrast with the heat equation for which non-negative solutions are automatically positive everywhere on its domain of definition, and is an argument of the more physical soundness of PME for diffusion process modelization.

Since Eq.~\eqref{eq:PME} is not defined for $u<0$, it is typically reformulated for numerical solving as follows:
\begin{equation}
    \dfrac{\partial u}{\partial t} = \nabla \cdot \left(\kappa |u|^m \nabla u\right), \quad m > 0\\
  \label{eq:posPME}
\end{equation}
called the signed PME. Solving Eq.~\eqref{eq:posPME} while ensuring non-negativity of the solution over the whole domain $\Omega$ is equivalent to solving  Eq.~\eqref{eq:PME}. Conventional numerical methods, such as finite element methods (FEM), face notable challenges when applied to Eq.~\eqref{eq:posPME}. With $\vec n$ denoting the normal vector at the interface $\Gamma(u)$, the quantity $\nabla u \cdot {\vec n}$ may lack regularity, depending on $m$ and $u_0$. This irregularity in the solution's derivatives implies that failure to capture $\Gamma(u)$ sharply will induce spurious oscillations, compromising non-negativity and accuracy.

To overcome these challenges, we propose here the use of \xm{} method which is specifically designed to handle problems with sharp interfaces in the FEM context. \xm{} has demonstrated its effectiveness in scenarios such as the Stefan problem or two-phase flows \cite{moes2023extreme, quiriny2024x}, particularly in accurately capturing complex interface dynamics without the need for enriched elements or mesh topology modifications. We will show that \xm{} is a valuable tool for accurately simulating the PME while addressing some limitations of traditional FEM approaches.

\section{The Porous Medium Equation}
\label{sec:PME}

\begin{figure}[h!]
  \begin{center}
    \includegraphics[width=0.48\textwidth]{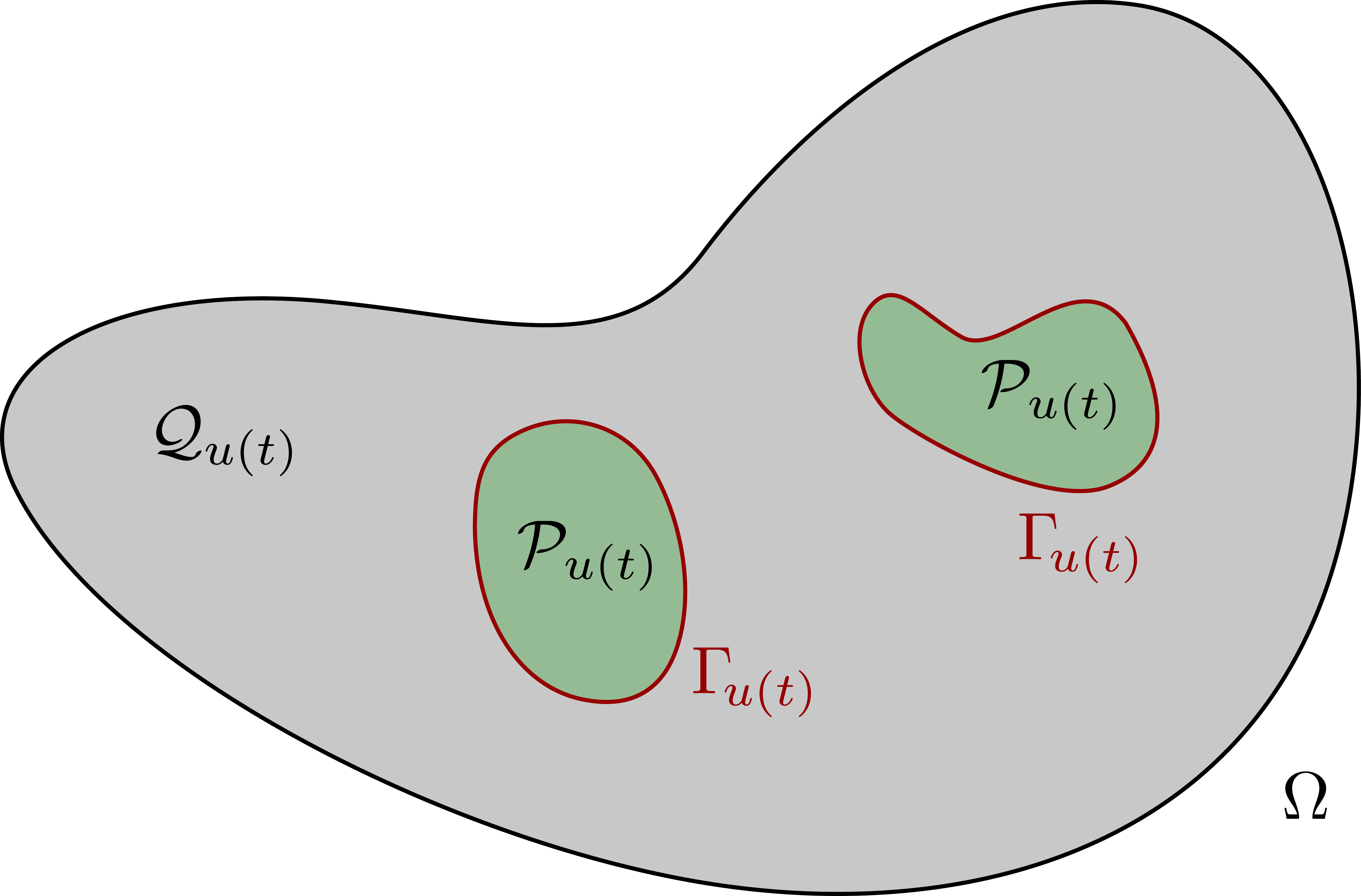}
  \end{center}
  \caption{Green areas correspond to \emph{the phase} $\mathcal{P}_{u(t)}$ (part of the domain where $u(t) > 0$). The grey area corresponds to the \emph{empty region} $\mathcal{Q}_{u(t)}$ (part of the domain where $u(t) = 0$) and red lines to \emph{the interface} $\Gamma(u(t))$.}
  \label{fig:defPQ}
\end{figure}

In this section, we will detail important properties of PME and state-of-the-art methods to solve it. We consider a domain $\Omega \subset \mathbb{R}^d$ and a known function $u_{0}\in\mathcal{C}^0(\Omega)$ such that:
\begin{equation}
  \left\{
  \begin{array}{lcll}
    \vspace{0.5em}
    u_{0}(\vx) & \geq & 0, & \forall \vx\in\Omega\\
    \nabla u_0(\vx) \cdot \vec{n} & = & 0, & \forall \vx\in\partial\Omega\\
  \end{array}
  \right.
  \label{eq:initSol}
\end{equation}
The objective is, for fixed $m>0$ and $\kappa>0$, to find $u$ such that:
\begin{equation}
  \left\{
  \begin{array}{lcll}
    \vspace{0.5em}
    \dpfrac{u}{t}(\vec{x}, t) & = & \nabla \cdot \left(\kappa |u(\vx, t)|^m \nabla u(\vx, t)\right) & ,  ~~(\vx, t)\in\Omega\times [t_0, T]\\
    \vspace{0.5em}
    u(\vx, t_0) & = & u_0(\vx) & , ~~\vx\in\Omega \\
    \nabla u(\vx,t) \cdot \vec{n} & = & 0 & ,  ~~(\vx, t)\in\partial\Omega\times [t_0, T] \\
    u(\vx, t) & \geq & 0 & ,  ~~(\vx, t)\in\Omega\times [t_0, T] \\
  \end{array}
  \right.
  \label{eq:posPME2}
\end{equation}
As mentioned in the introduction, solutions of the PME can exhibit an interface separating regions where the solution is null from region where the solution is strictly positive. Therefore, problems modeled with the PME fall into the category of \onep\, problems. For clarity purposes, we are defining in this paper $\mathcal{P}_{u(t)} = \{ \vx\in\Omega \text{ such as } u(\vx, t) > 0\}$ , $\Gamma(u(t)) = \partial\mathcal{P}_{u(t)}$ and $\mathcal{Q}_{u(t)} = \Omega\setminus (\mathcal{P}_{u(t)}\cup\partial\mathcal{P}_{u(t)})$. $\mathcal{P}_{u(t)}$ will be refered to as \emph{the phase}, $\Gamma(u(t))$ as \emph{the interface} and $\mathcal{Q}_{u(t)}$ the \emph{empty region} for a given solution $u(t)$ on $\Omega$ (Figure~\ref{fig:defPQ}).

PME presents various notable properties that a numerical simulation must be able to ensure. We will highlight the most crucial ones, and a more extensive list can be found in \cite{vazquez2007porous}.

\subsection{Mass Conservation}

An important property of PME is the invariance in time of the quantity $\int_\Omega u \,d\Omega$ when homogeneous Neumann boundary conditions are applied on $\partial\Omega$. A comprehensive proof is available in \cite{vazquez2007porous}. In the context of PME modeling for non-linear heat transfer, this corresponds to the conservation of internal energy. In the case of PME modeling for gas flow in a porous medium, it signifies mass conservation. For convenience, we refer to this property as mass conservation in this paper.

\subsection{Interface Velocity}

Using Reynolds transport theorem and the mass conservation property of PME's solutions, it is possible to derive the interface normal velocity $\vec{v}_{\iGamma}$ \cite{vazquez2007porous}:

\begin{equation}
  \vec{v}_{\iGamma(u(t))}(\vec{x}, t) = \lim_{\vec{x} \in \mathcal{P}_{u(t)} \rightarrow \Gamma(u(t))} -\kappa \nabla \left( \frac{u^m(\vec{x}, t)}{m}\right)
  \label{eq:intVelocity}
\end{equation}

This equation highlights that there exists a set of initial conditions for which the interface $\Gamma(u(t))$ will not move for a finite amount of time. Indeed, choosing $u_0$ such that $\nabla (u^m_0)$ is continuous and vanishes on $\Gamma(u_0)$ leads to $\vec{v}_{\iGamma(u_0)}(\vec{x}, t_0) = \vec{0}$. This phenomenon is known as the \emph{waiting time} property in the literature. A robust numerical simulation of the PME should accurately reproduce this specific behavior.

\subsection{Barenblatt-Pattle solution}

Another noteworthy result is the existence of analytical solutions of the form:
\begin{equation}
  u_{BP}(\vec{x}, t) = t^{-\alpha}\left((C-k\|\vec{x}\|^2 t^{-2\beta})^+\right)^{\frac{1}{m}}
  \label{eq:barenblatt}
\end{equation}
with
\begin{equation}
  \left\{
  \begin{array}{lll}
    \vspace{0.5em}
    \alpha & = & \dfrac{d}{md+2}\\
    \vspace{0.5em}
    \beta & = & \dfrac{\alpha}{d}\\
    k & = & \dfrac{m\alpha}{2\kappa d}\\
  \end{array}
  \right.
\end{equation}
and $C$ an arbitrary constant. This solution correspond to an initial condition $u_0$ being a Dirac distribution centered on $\vec 0$ for which the total mass $\int_\Omega u_0 \id\Omega = M(C)$ depends on $C$. This analytical solution will be used to characterize the convergence rate of the proposed numerical method. The interface $\Gamma(u_{BP}(t))$ is well-defined as:
\begin{equation}
  \Gamma(u_{BP}(t)) = \left\{ \vec x\in\Omega \mid \|\vec{x}\| = f(t), f(t) = \left(\dfrac{C}{k}\right)^{\frac{1}{2}}t^\beta\right\}
  \label{eq:bpIntF}
\end{equation}

\subsection{Solving PME with Finite Element Method}

\begin{figure}[h!]
  \begin{center}
    \includegraphics[width=0.48\textwidth]{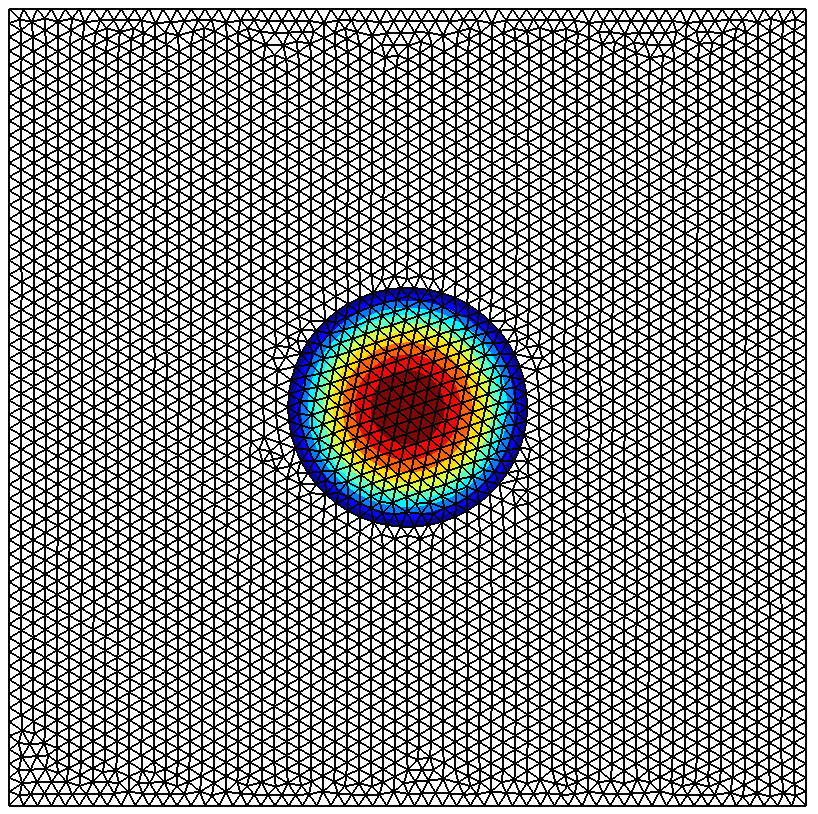}
    \hfill
    \includegraphics[width=0.48\textwidth]{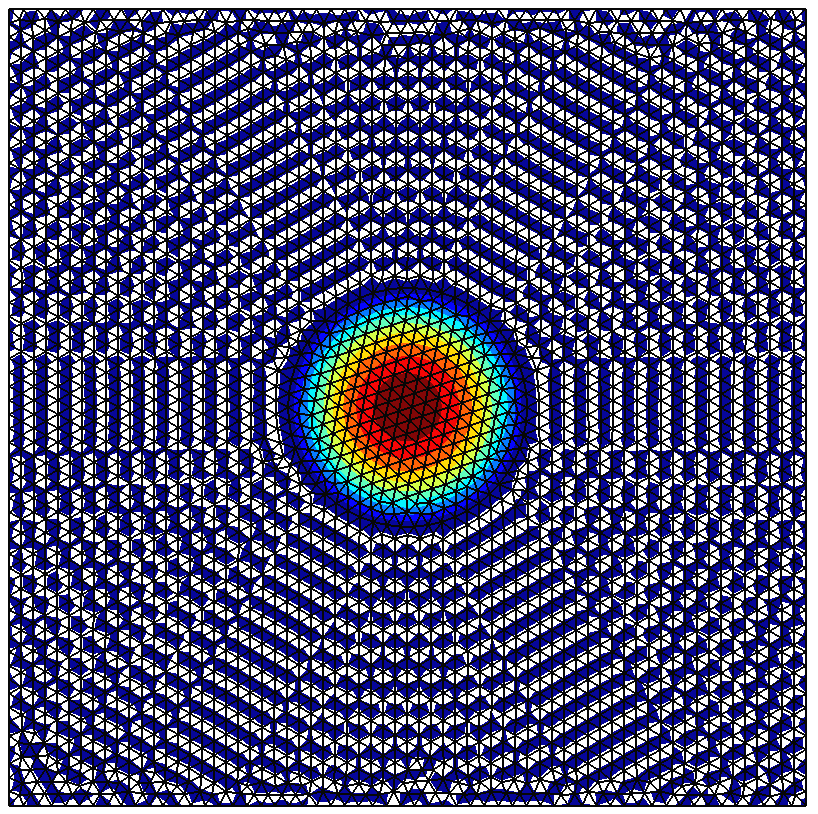}
  \end{center}
  \caption{Numerical solution of the PME obtained with a classical FEM. Colored areas correspond to the part of the domain $\Omega$ where $u\geq 0$. Left: initial condition at $t=t_0$. Right: numerical solution at $t_1 = t_0 + \Delta t$.}
  \label{fig:exOsc}
\end{figure}

When solving Eq.~\eqref{eq:posPME2} with FEM, the property of non-negativity is violated, and oscillations appear at interface $\Gamma(u(t))$, propagating throughout the entire domain $\Omega$. An example of such oscillations is illustrated Figure~\ref{fig:exOsc} for $d=2$.

The manifestation of this behavior coupled with a loss of regularity in the vicinity of the interface, results in a diminished convergence rate of the numerical scheme. The theorical convergence rate using a uniform linear finite element discetization of size $h$ in space and a backward-Euler discretization in time established in \cite{ebmeyer1998error} is:

\begin{equation}
  \|u^h - u\|_{L^2([t_0, T];L^2(\Omega))} \leq Ch^{\frac{m^2+6m+8}{6m^2+14m+8}}, \text{ with } \Delta t = \mathcal{O}(h^{\frac{5m+4}{2m}})
  \label{eq:convRate}
\end{equation}
As $m$ tends to $0$, the convergence rate reaches its maximum of $1$ and gradually approaches $\frac{1}{6}$ with increasing values of $m$.

To circumvent this issue, several approaches have been proposed. A first idea is to use a numerical scheme which ensure non-negativity of the solution \cite{liu2011high, zhang2009numerical, chen2016third}. For example, \cite{zhang2009numerical} proposes a method where a local discontinuous Galerkin finite element method is coupled with a custom made limiter to ensure non-negativity of the solution. A demonstration is provided showing that this limiter ensure non-negativity of the solution for the discontinuous $\mathbb{P}_0$ finite elements, and provide numerical 1D results for discontinuous $\mathbb{P}_2$ finite elements. However, these methods are adding numerical viscosity leading to a loss of precision of the numerical solution near the interface.

A second type of approach is to use a time-dependent spatial discretization. In \cite{ngo2017study} is proposed an adaptative moving mesh finite element method. Meshes used for solving the PME are time-dependent but all have the same topology (same number of vertices and elements, and same connectivity). Each mesh is generated such as its elements sizes respect an imposed metric and different metric are investigated. They are showing that using a target metric based on the Hessian of the solution at previous time step allows to retreive a convergence rate in $\mathcal{O}(h^2)$ for $m=1$ and $m=2$ for the Barenblatt-Pattle solution. It is important to note that with this method, non-negativity of the solution is broken and oscillations at the interface $\Gamma(u(t))$ are still showing. However, the target metric used leads to a smaller element size near the interface decreasing significantly oscillations amplitude and propagation in the whole computational domain.

In \cite{baines2005moving} is proposed a moving mesh finite element algorithm for PDE with moving boundaries, and is applied to the Barenblatt-Pattle solution of the PME. In that case, only the part of the computational domain where $u>0$ is meshed and the mesh evolves in time in order to ensure mass conservation of the solution. This allows to obtain a second-order convergence rate for $m=1$ and a first-order convergence rate for $m=3$. However, this approach does not handle topology changes which can lead the mesh to overlap itself and need to be coupled with a remeshing algorithm for such cases.

To the best of our knowledge, all numerical methods for solving the PME described in the literature involve violating either the non-negativity of the solution, mass conservation, or modifying the equation through direct or indirect regularization.

The methodology presented in this paper falls within the moving meshes category. The principle of the \xm{} method is to determine, at each time step $t_n$, the solution $u_{n}$ and the spatial discretization $\mathcal{T}_n$ in such a way that $\mathcal{T}_n$ maintains the same topology as $\mathcal{T}_{n-1}$ and fully represents $\Gamma(u_n)$ through its vertices and edges. Key criteria for evaluating the method's effectiveness include ensuring non-negativity at each time step, improving the convergence rate with respect to the spatial mesh size $h$, enhancing precision in interface localization, maintaining robustness in the face of topology changes in $\mathcal{P}_u(t)$, and preserving the mass conservation property of the PME.

Firstly, the \xm{} approach will be presented in more detail, followed by a comprehensive description of the moving mesh strategy. Finally, the results and relevance of the proposed method will be discussed.

\section{The \xm{} approach}
\label{sec:xmesh}

As previously mentioned, the goal of the \xm{} approach is to use a time-dependent spatial discretization. In order to develop a discrete formulation of the problem with this constraint, we introduce the arbitrary Lagrangian Eulerian frame of reference (\cite{boffi2004stability, donea2004arbitrary}). We define the reference domain $\Omega_0$ as $\Omega_0 = \Omega$, and a mapping $\vec{X} \in H^1(\Omega_0\times[t_0, T])$:
\begin{equation}
  \begin{array}{rll}
    \vec{X} : \Omega_0\times [t_0, T] & \rightarrow & \mathbb{R}^d\\
    (\vec{X}_0, t) & \mapsto & \vec X(\vec{X}_0, t)
  \end{array}
\end{equation}
We define $\Omega_t = \{ \vec{X}(\vec{X}_0, t), \vec{X}_0\in\Omega_0, t\in[t_0, T]\}$, $\vec{w} = \left. \dpfrac{\vec{X}}{t}\right|_{\vec{X}_0}$, $\vec{F} = \left. \dpfrac{\vec{X}}{\vec{X}_0}\right|_t$ and $J = \text{det }{\vec{F}}$.
Assumptions on the mapping are the following:
\begin{equation}
  \left\{
  \begin{array}{l}
    \vspace{0.5em}
    \mathbf{X}(\vec{X}_0, t_0) = \vec{X}_0\text{,}\\
    \vspace{0.5em}
    \text{the domain shape is unchanged, that is to say } \Omega_t = \Omega_0\text{, }\forall t\in[t_0, T],\\
    J(\vec{X}_0, t) \geq 0 \text{, } \forall (\vec{X}_0, t)\in\Omega_0\times[t_0, T] \text{.}\\
  \end{array}
  \right.
  \label{eq:condX}
\end{equation}
To obtain a weak formulation of Eq.~\eqref{eq:posPME2}, we define the test function space \( \ts{U} \) as follows, assuming functions within it have proper regularity on \( \Omega \times [t_0, T] \):
\begin{equation}
  \ts{U} = \{ \tf u : \Omega \times [t_0, T] \rightarrow \mathbb{R} \mid \left. \dpfrac{\tf u}{t} \right|_{\vec{X}} = 0 \}
  \label{eq:defTestSpace}
\end{equation}
With $\bar{u}\in\ts{U}$, Eq.~\eqref{eq:posPME2} gives:
\begin{eqnarray}
  \begin{array}{crcl}
    & \dfrac{\text{d}}{\text{d}t} \dint_\Omega u\bar{u} \id\Omega & = & \dint_\Omega \dfrac{\partial u}{\partial t}\bar{u} \id\Omega - \dint_\Omega u\vec{w}.\nabla \bar{u}\id\Omega = \dint_\Omega \nabla.(\kappa|u|^m\nabla u)\bar{u} \id\Omega - \dint_\Omega u\vec{w}\nabla \bar{u}\id\Omega\\
    \\
    \Leftrightarrow & \dfrac{\text{d}}{\text{d}t} \dint_\Omega u\bar{u} \id\Omega & = & - \dint_\Omega \kappa|u|^m \nabla u.\nabla \bar{u} \id\Omega - \dint_\Omega u\vec{w}\nabla \bar{u}\id\Omega.\\
  \end{array}
  \label{eq:varform}
\end{eqnarray}
The time integration of Eq.~\eqref{eq:varform} between two instants $t_n$ and $t_{n+1}$ is performed using a $\theta$ scheme. Noting $\Delta t_n = t_{n+1} - t_n$ and $u_n(\vec{X}) = u(\vec{X}, t_n)$ we have:

\begin{equation}
  \dint_\Omega u_{n+1}\bar{u} \id\Omega - \dint_\Omega u_{n}\bar{u} \id\Omega \approx \Delta t_n \left(\theta\dfrac{\text{d}}{\text{d}t} \dint_\Omega u_{n+1}\bar{u} \id\Omega + (1-\theta)\dfrac{\text{d}}{\text{d}t} \dint_\Omega u_n\bar{u} \id\Omega\right) \text{, } \theta\in[0, 1]
  \label{eq:timediscr}
\end{equation}
In this paper, unless mentioned otherwise, $\theta$ is set to $\frac{1}{2}$.

The spatial discretization is done using an order $1$ finite element method. The reference domain $\Omega_0$ is partionned into a triangular mesh $\mathcal{T}_0$ for which the set of nodes is denoted $\mathcal{N}_0$. The spatial discretization at $t$ is the triangulation $\mathcal{T}(t)$ obtained from applying the mapping $\vec X$ to $\mathcal{T}_0$. By construction, the quantity $\mathbf{w}$ corresponds to the mesh velocity and the mesh topology (number of nodes and connectivities) does not depend on time. The set of nodes of $\mathcal{T}(t)$ is denoted $\mathcal{N}$, and we note $\vec{x}^p(t)$ the coordinates of $p\in\mathcal{N}$. For convenience purposes, we will also define:
\begin{equation}
  \vec{X}(p, t) \equiv \vec{x}^p(t) = \vec{X}(\vec{x}^p(0), t)\text{, }p\in\mathcal{N}
  \label{eq:defXp}
\end{equation}
The discrete solution $u^h$ at instant $t_n$ is noted $u^h_n \equiv u^h_{\mid_{t_n}}$ and verifies:
\begin{equation}
  u^h_n(\vec{X}) = \sum_{i\in\mathcal{N}} U^i_n \phi_i(\vec{X})
  \label{eq:spacediscr}
\end{equation}
where $\phi_i$ is the finite element linear basis function associated to node $i$ and $U^i_n$ the value of the discrete scalar field $u^h_n$ at node $i$.
The test function $\bar{u}$ is discretized using the same process. The discrete test function $\bar{u}^h$ is:
\begin{equation}
  \bar{u}^h(\vec{X}) = \sum_{i\in\mathcal{N}} \bar{U}^i \phi_i(\vec{X}).
  \label{eq:spacediscrtest}
\end{equation}
and the discrete test function space $\ts{U}^h$ is:
\begin{equation}
  \ts{U}^h = \{ u:\Omega\times[t_0, T]\rightarrow\mathbb{R} \mid u(\vec{X}) = \sum_{i\in\mathcal{N}} U^i\phi_i(\vec{X}), U^i\in\mathbb{R}\text{ }\forall i\in\mathcal{N}\}
  \label{eq:defDiscrTestSpace}
\end{equation}

With the spatial discretization described above and the time integration method outlined in Eq.~\eqref{eq:timediscr} over the time interval $[t_n,t_{n+1}]$, Eq.~\eqref{eq:varform} becomes:
\begin{eqnarray}
  \begin{array}{lll}
    r(\bar{u}^h, u^h_{n+1},\vec{X}_{n+1}, u^h_{n},\vec{X}_{n}) & \equiv & \dint_\Omega u^h_{n+1}\bar{u}^h \id\Omega - \dint_\Omega u^h_{n}\bar{u}^h \id\Omega\\
    & +& \Delta t_n \theta \left(\dint_\Omega \kappa|u^h_{n+1}|^m \nabla u^h_{n+1}.\nabla \bar{u}^h \id\Omega + \dint_\Omega u^h_{n+1}\vec{w}_{n+1}\nabla \bar{u}^h\id\Omega \right) \\
    & +& \Delta t_n (1-\theta) \left(\dint_\Omega \kappa|u^h_{n}|^m \nabla u^h_{n}.\nabla \bar{u}^h \id\Omega + \dint_\Omega u^h_{n}\vec{w}_{n}\nabla \bar{u}^h\id\Omega \right) \\
    & = & 0 \text{, } \forall \bar{u}^h\in\ts{U}^h\\
  \end{array}
  \label{eq:defr}
\end{eqnarray}

The \xm{} approach consists in, knowing $(u^h_{n}, \vec{X}_n)$, finding $(u^h_{n+1}, \vec{X}_{n+1})$ such as:
\begin{equation}
  \left\{
  \begin{array}{l}
    \vspace{0.5em}
    \vec{X}_{n+1} \text{ respects conditions Eq.~\eqref{eq:condX}}\\
        
    \vspace{0.5em}
    u^h_{n+1}(\vec{X}_{n+1}) \ge 0 \text{, } \forall (\vec{X}_{n+1})\in\Omega\\
    
    r(\bar{u}^h, u^h_{n+1},\vec{X}_{n+1}, u^h_{n},\vec{X}_{n}) = 0 \text{, } \forall \bar{u}^h\in\ts{U}^h\\
  \end{array}
  \right.
\end{equation}
It is important to note that this set of conditions is equivalent to:
\begin{equation}
  \left\{
  \begin{array}{l}
    \vspace{0.5em}
    \vec{X}_{n+1} \text{ respects conditions Eq.~\eqref{eq:condX}}\\
        
    \vspace{0.5em}
    \Gamma_{n+1}^h \equiv \Gamma(u^h_{n+1}) \text{ is fully represented with vertices of } \mathcal{T}_{n+1}\\
    
    r(\bar{u}^h, u^h_{n+1},\vec{X}_{n+1}, u^h_{n},\vec{X}_{n}) = 0 \text{, } \forall \bar{u}^h\in\ts{U}^h\\
  \end{array}
  \right.
  \label{eq:pbXM0}
\end{equation}
which can be rewritten as:
\begin{equation}
  \left\{
  \begin{array}{l}
    \vspace{0.5em}
    \vec{X}_{n+1} \text{ respects conditions Eq.~\eqref{eq:condX}}\\
        
    \vspace{0.5em}
    \Gamma_{n+1}^h \text{ is fully represented with vertices of } \mathcal{T}_{n+1}\\
    
    R((U_{n}, \vec{X}_{n}), (U_{n+1}, \vec{X}_{n+1})) = \left(\displaystyle\sum_{i\in\mathcal{N}}r^2(\phi_i, u^h_{n+1},\vec{X}_{n+1}, u^h_{n},\vec{X}_{n})\right)^{\frac{1}{2}} = 0\\
  \end{array}
  \right.
  \label{eq:pbXM}
\end{equation}
In the following sections will be detailed the general method to solve problem \eqref{eq:pbXM}.

\subsection{Solving scheme}

\begin{figure}[h!]
  \begin{center}
    \includegraphics[width=0.95\textwidth]{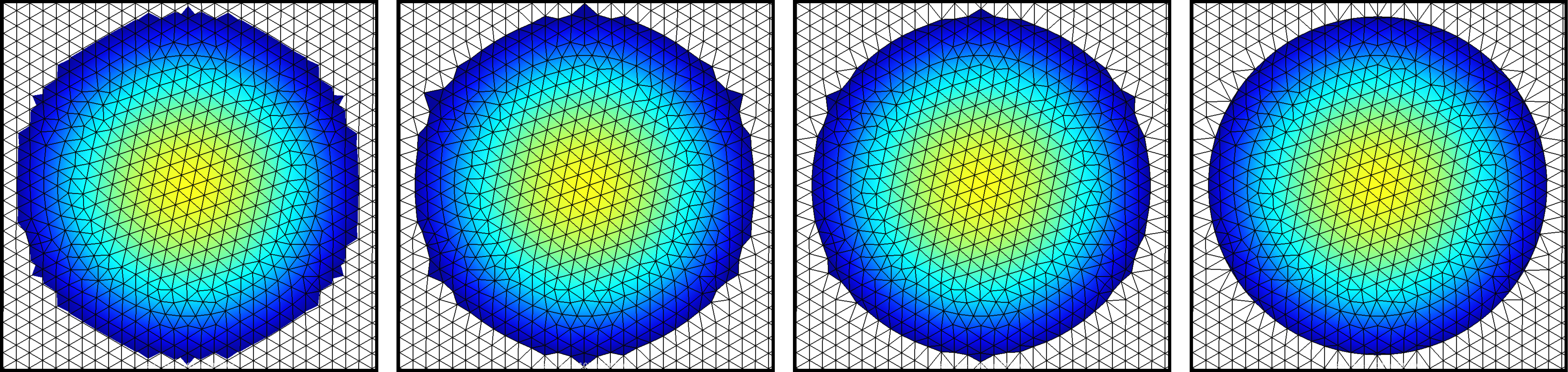}
  \end{center}
  \caption{Evolution of $(U_{n+1}^k, \vec X_{n+1}^k)$ when applying Algorithm \ref{alg:xm}. Colored areas correspond to part of the domain where $u^h_{n+1}>0$. On the left is represented the solution $U^0_{n+1}$ obtained on the mesh intial guess $\vec X_{n+1}^0 = \vec X_0$. Middle left and right present the mesh $\vec X^k_{n+1}$ and associated solution $U^k_{n+1}$ for $k=1$ and $k=2$ respectively. On the right is depicted the converged solution $(U_{n+1}, \vec X_{n+1})$.}
  \label{fig:exUpdates}
\end{figure}

We assume that the valid solution $(u^h_n, \vec{X}_n)$ at time step $t_n$ is known. The method developed to solve the non-linear problem described in Eq.~\eqref{eq:pbXM} is outlined in Algorithm \ref{alg:xm}. The general approach begins by selecting nodes $\mathcal{N}(\Gamma^h_{n+1})$, which will be used to define the interface $\Gamma^h_{n+1}$ at time step $t_{n+1}$. Next, the pair $(U_{n+1}, \vec X_{n+1})$ that satisfies Eq.~\eqref{eq:pbXM} is determined through successive updates $k$ of the mesh $\vec X^k_{n+1}$ and the corresponding solution $U^k_{n+1}$. This iterative process continues until the residual $R((U_n, \vec X_n), (U^k_{n+1}, \vec X^k_{n+1}))$ meets a specified tolerance.
\begin{algorithm}[h!]
  \caption{Solving scheme for time step $t_{n+1}$. Solution $(U_n, \vec{X}_n)$ at instant $t_n$ is known. $tol$ is a user-imposed convergence tolerance.}
  $k \gets 0$\;
  $\vec X^0_{n+1} \gets \vec X_{0}$\;
  $U_{n+1}^{0} \gets solve((U_n, \vec{X}_n), \vec{X}_{n+1}^{0})$\;
  $\mathcal{N}(\Gamma^h_{n+1}) \gets localize\_interface(U_{n+1}^{0}, \vec{X}_{n+1}^{0})$\;
  \While{$\epsilon > tol$}{
    $\vec{X}_{n+1}^{k+1} \gets update\_interface(\mathcal{N}(\Gamma^h_{n+1}), (U_n, \vec{X}_n) , (U_{n+1}^{k}, \vec{X}_{n+1}^{k}))$\;
    $U_{n+1}^{k+1} \gets solve((U_n, \vec{X}_n), \vec{X}_{n+1}^{k+1})$\;
    $\epsilon \gets R((U_{n}, \vec{X}_{n}), (U_{n+1}^{k+1}, \vec{X}_{n+1}^{k+1}))$\;
    $k \gets k + 1$\;
  }
  $(U_{n+1}, \vec{X}_{n+1}) \gets (U_{n+1}^{k}, \vec{X}_{n+1}^{k})$\;
  \label{alg:xm}
\end{algorithm}
Figure~\ref{fig:exUpdates} present the evolution of $(U_{n+1}^k, \vec X_{n+1}^k)$ and the converged result obtained when applying Algorithm \ref{alg:xm}. In this example, expected interface $\Gamma^h_{n+1}$ is circular and cannot be represented with mesh edges from the initial guess $\vec X^0_{n+1}$. Successive mesh updates converge to a mesh $\vec X_{n+1}$ from which a better representation of $\Gamma^h_{n+1}$ is obtained using only mesh edges.

In order for this approach to be effective, the $solve$ method has to provide a solution from which selecting nodes to belong to $\Gamma^h_{n+1}$ is straightforward (i.e.\ a solution without oscillations propagating in the whole domain). We also need to be able to estimate the interface $\Gamma^h_{n+1}$ location despite not having an explicit description of it, and finally to properly define the mesh velocity $\vec w$.

\subsection{Computation of solution $U^k_{n+1}$ associated to $\vec{X}_{n+1}^k$}

The step referred to as $solve$ in Algorithm \ref{alg:xm} is detailed here. This step aims to compute $U_{n+1}^k$ for a fixed $\vec{X}_{n+1}^k$ when $(U_n, \vec{X}_{n})$ is known. To streamline the notations, in this section $U_{n+1}^k$ and $\vec{X}_{n+1}^k$ will be respectively noted $U$ and $\vec X$, and $\phi_i$ will refer to the first order finite element basis function associated to node $i\in\mathcal{N}$.

This problem is non-linear and can, for example, be solved using a Newton-Raphson scheme. An initial guess $U_0$ is made and the solution is updated following:
\begin{equation}
  H_j.U_{j+1} = H_j.U_{j} - G_j
  \label{eq:NR}
\end{equation}
while $\|G_j\|$ is greater than a user imposed tolerance, with:
\begin{equation}
  \left\{
  \begin{array}{lcl}
    \vspace{0.5em}
    G_j & = &
  \left[
  \begin{array}{c}
    r(\phi_1, U_{j}, \vec X, U_n, \vec X_n)\\
    \vdots\\
    r(\phi_N, U_{j}, \vec X, U_n, \vec X_n)\\
  \end{array}
  \right]
    \\
    H_j & = & \dpfrac{G_j}{U_j}\\
  \end{array}
  \right.
\end{equation}
As mentioned in Section~\ref{sec:PME}, if nodes of $\mathcal{T}_{n+1}$ are not describing $\Gamma^h_{n+1}$, solving Eq.~\eqref{eq:NR} will lead to a solution that violates non-negativity and exhibits oscillations propagating throughout the entire domain.

\begin{figure}[h!]
  \begin{center}
    \includegraphics[width=0.48\textwidth]{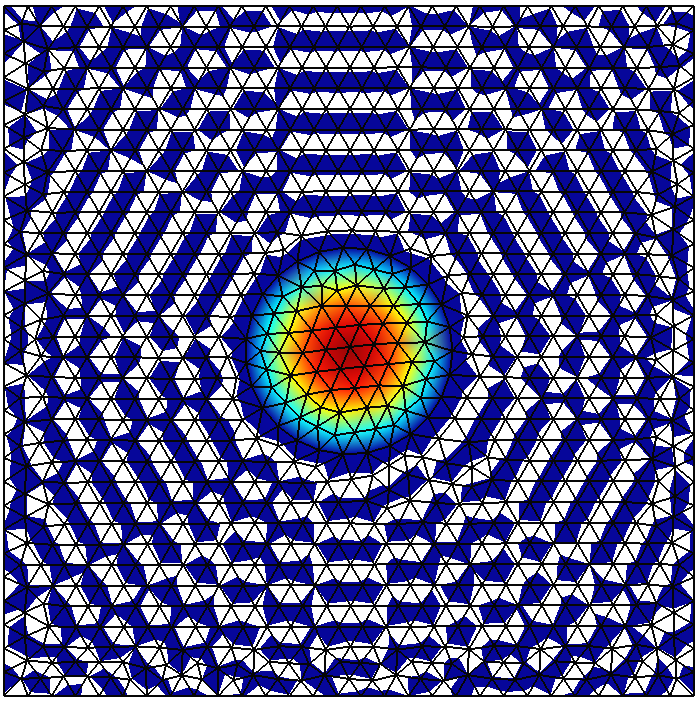}
    \hfill
    \includegraphics[width=0.48\textwidth]{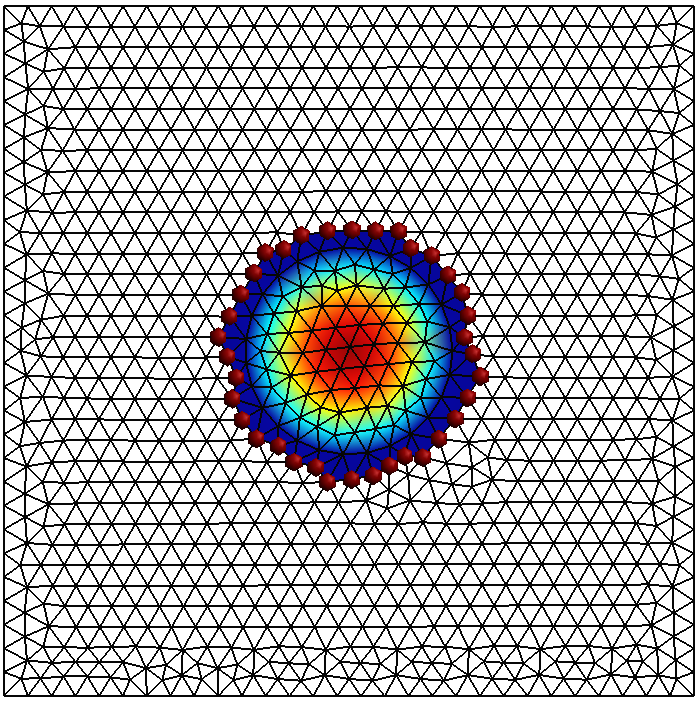}
  \end{center}
  \caption{Numerical solutions obtained at $t_{n+1}$ when $\mathcal{T}_{n+1}$ does not fully describe $\Gamma^h_{n+1}$. Colored areas correspond to the part of the domain $\Omega$ where $u > 0$. Left: solution obtained with classical Newton-Raphson iterations. Right: solution obtained with Newton-Raphson iterations with Lagrange multipliers described Eq.~\eqref{eq:NRmod}. The red spheres highlight nodes associated to a Lagrange multiplier.}
  \label{fig:NR}
\end{figure}

In order to circumvent this issue and ensure non-negativity of the converged solution $U$, the Newton-Raphson iteration Eq.~\eqref{eq:NR} is modified by introducing Lagrange multipliers $\lambda^\mathcal{D}$ to ensure an homogeneous Dirichlet boundary condition on a well chosen set of vertices $\mathcal{D}\subset\mathcal{N}$. Defining $\mathcal{F} = \mathcal{N}\setminus \mathcal{D}$ and $I$ the identity matrix, the Newton-Raphson iteration with Lagrange multipliers is written:
\begin{equation}
  \left[
  \begin{array}{ccc}
    H_j^{\mathcal{F}\mathcal{F}} & H_j^{\mathcal{F}\mathcal{D}} & 0 \\
    H_j^{\mathcal{D} \mathcal{F}} & H_j^{\mathcal{D}\mathcal{D}} & -I \\
    0 & -I & 0 \\
  \end{array}
  \right]
  \left[
  \begin{array}{c}
    U_{j+1}^{\mathcal{F}}\\
    U_{j+1}^{\mathcal{D}}\\
    \lambda^\mathcal{D}_{j+1}\\
  \end{array}
  \right]
  =
  \left[
  \begin{array}{ccc}
    H_j^{\mathcal{F}\mathcal{F}} & H_j^{\mathcal{F}\mathcal{D}} & 0 \\
    H_j^{\mathcal{D} \mathcal{F}} & H_j^{\mathcal{D}\mathcal{D}} & 0 \\
    0 & 0 & 0 \\
  \end{array}
  \right]  
  \left[
  \begin{array}{c}
    U_j^{\mathcal{F}}\\
    U_j^{\mathcal{D}}\\
    0\\
  \end{array}
  \right]
  -
  \left[
  \begin{array}{c}
    G_j^{\mathcal{F}}\\
    G_j^{\mathcal{D}}\\
    0\\
  \end{array}
  \right]
  \label{eq:NRmod}
\end{equation}
It should be noted that the application of Lagrange multipliers is targeted solely at nodes where the solution $U$ would be negative, necessitating the retention of strictly positive multipliers. Thus, to compute the solution update $U_{j+1}$, we define the set of vertices $\mathcal{D}$ in the following manner:
\begin{equation}
  \mathcal D = \{m\in\mathcal{N}\mid U_j^m < 0 \text{ or } \lambda^m_j>0\}
\end{equation}
As we know that $U_{j+1}^{\mathcal{D}} = 0$, the system of equations \ref{eq:NRmod} can be rewritten in a more convenient way:
\begin{equation}
  \left[
  \begin{array}{cc}
    H_j^{\mathcal{F}\mathcal{F}}  & 0 \\
    H_j^{\mathcal{D} \mathcal{F}} & -I \\
  \end{array}
  \right]
  \left[
  \begin{array}{c}
    U_{j+1}^{\mathcal{F}}\\
    \lambda_{j+1}\\
  \end{array}
  \right]
  =
  \left[
  \begin{array}{cc}
    H_j^{\mathcal{F}\mathcal{F}} & H_j^{\mathcal{F}\mathcal{D}}\\
    H_j^{\mathcal{D} \mathcal{F}} & H_j^{\mathcal{D}\mathcal{D}}\\
  \end{array}
  \right]  
  \left[
  \begin{array}{c}
    U_j^{\mathcal{F}}\\
    U_j^{\mathcal{D}}\\
  \end{array}
  \right]
  -
  \left[
  \begin{array}{c}
    G_j^{\mathcal{F}}\\
    G_j^{\mathcal{D}}\\
  \end{array}
  \right]    
  \label{eq:NRmod2}
\end{equation}
This resolution method, summarized in Algorithm \ref{alg:NRmod}, allows to ensure non-negativity of the solution everywhere and prevents oscillations to propagate in the whole domain $\Omega$, as highlighted Figure~\ref{fig:NR}.
\begin{algorithm}[h!]
  \caption{$solve$ method.}
  \textbf{Input:} $(U_n, \vec{X}_n) , \vec{X}$\;
  $j \gets 0$\;
  
  \While{$\epsilon > tol$}{
    \uIf{$j = 0$}{
      $\mathcal{D} = \emptyset$\;
    }
    \Else{
      $\mathcal{D} = \{m\in\mathcal{N}\mid U_j^m < 0 \text{ or } \lambda^m_j>0\}$\;
    }
    $\text{Compute } (U_{j+1}, \lambda^\mathcal{D}_{j+1}) \text{ verifying Eq.~\eqref{eq:NRmod2}} $\;
    $\epsilon = \left\|\dpfrac{}{\tf U} \left(r(\tf U, U_{j+1}, \vec X, U_n, \vec X_n) - (\tf{U}^\mathcal{D})^T.\lambda^\mathcal{D}_{j+1} \right)\right\|$\;
    $j\gets j+1$\;
  }
  $(U_{n+1}, \vec X_{n+1}) \gets (U_{n+1}^{k}, \vec X_{n+1}^{k})$\;
  \label{alg:NRmod}
\end{algorithm}

Introducing Lagrange multipliers modifies the problem from solving Eq.~\eqref{eq:defr} to addressing the following system instead:
\begin{equation}
  \left\{
  \begin{array}{lll}
    r(\bar{U}, U_{n+1}, \vec{X}_{n+1}, U_{n},\vec{X}_{n}) - (\tf{U}^\mathcal{D})^T.\lambda^\mathcal{D} & = & 0\text{, for all }\tf{U}\in\mathbb{R}^n\\
    U^D_{n+1} & = & \vec 0\\
  \end{array}
  \right.
  \label{eq:defrLag}
\end{equation}
Evaluating Eq.~\eqref{eq:defrLag} for the constant test function $\tf{u}^h = 1$, we obtain:
\begin{equation}
  \int_\Omega u^h_{n+1} \id\Omega - \int_\Omega u^h_{n}\id\Omega = \sum_{m\in\mathcal{D}} \lambda^m > 0.
  \label{eq:massConsLag}
\end{equation}
This highlights that, in addition to solving a different problem, mass conservation is no longer ensured.

However, it's essential to note that the solution obtained upon convergence of Algorithm \ref{alg:xm} will be non-negative on $\Omega$ without the use of Lagrange multipliers, resulting in $\mathcal{D} = \emptyset$.  Therefore, upon convergence of the \xm{} iterations, Eq.~\eqref{eq:defrLag} becomes equivalent to Eq.~\eqref{eq:defr}, ensuring mass conservation within the convergence tolerance set for Algorithm \ref{alg:xm}.

\subsection{Interface location estimation and mesh update}
This section focuses on estimating the location of the interface for a given solution $(U_{n+1}, \vec X_{n+1})$ and the subsequent mesh update required to capture interface $\Gamma^h_{n+1}$.

\begin{figure}[h!]
  \begin{center}
    \includegraphics[width=0.48\textwidth]{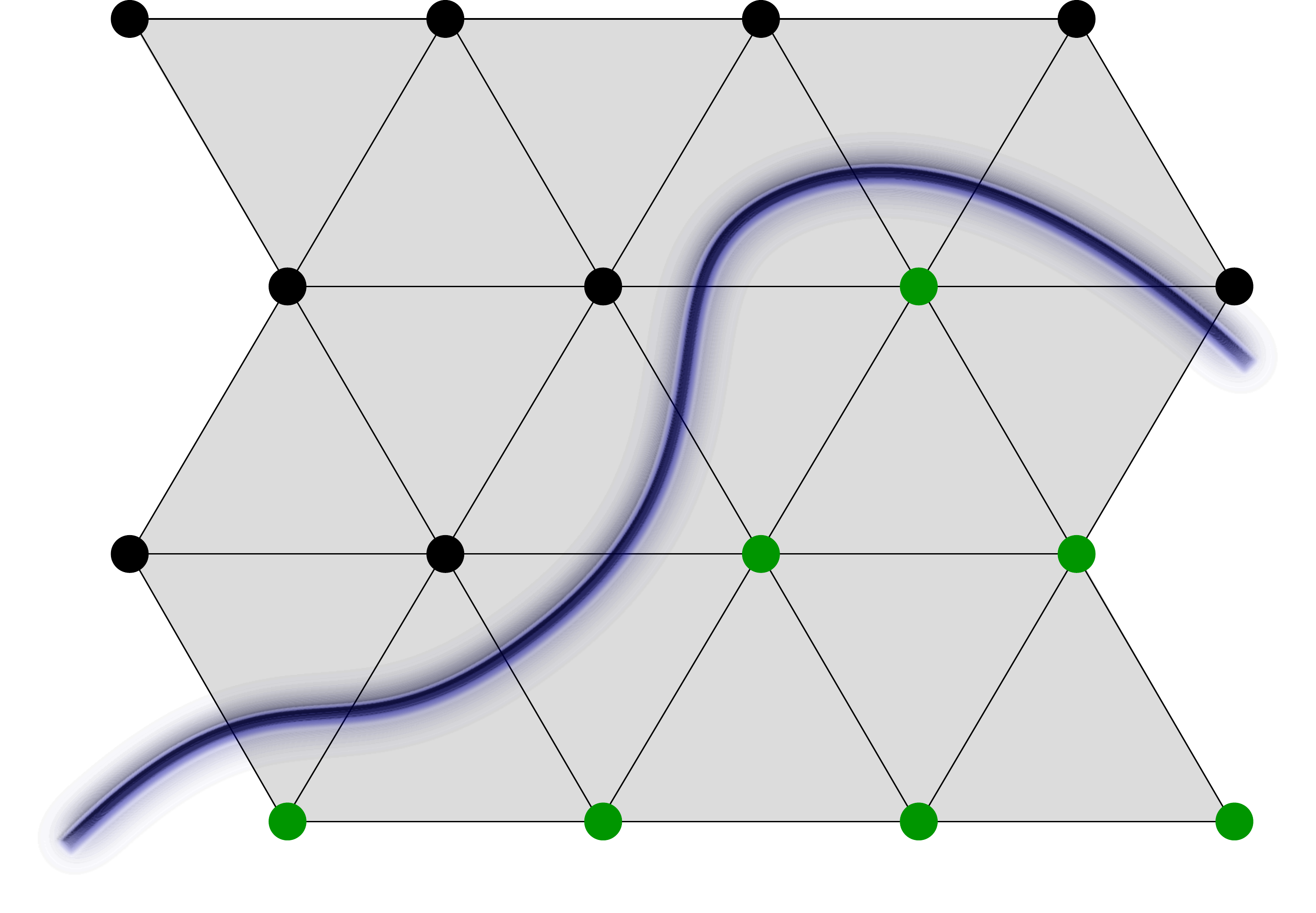}
    \hfill
    \includegraphics[width=0.48\textwidth]{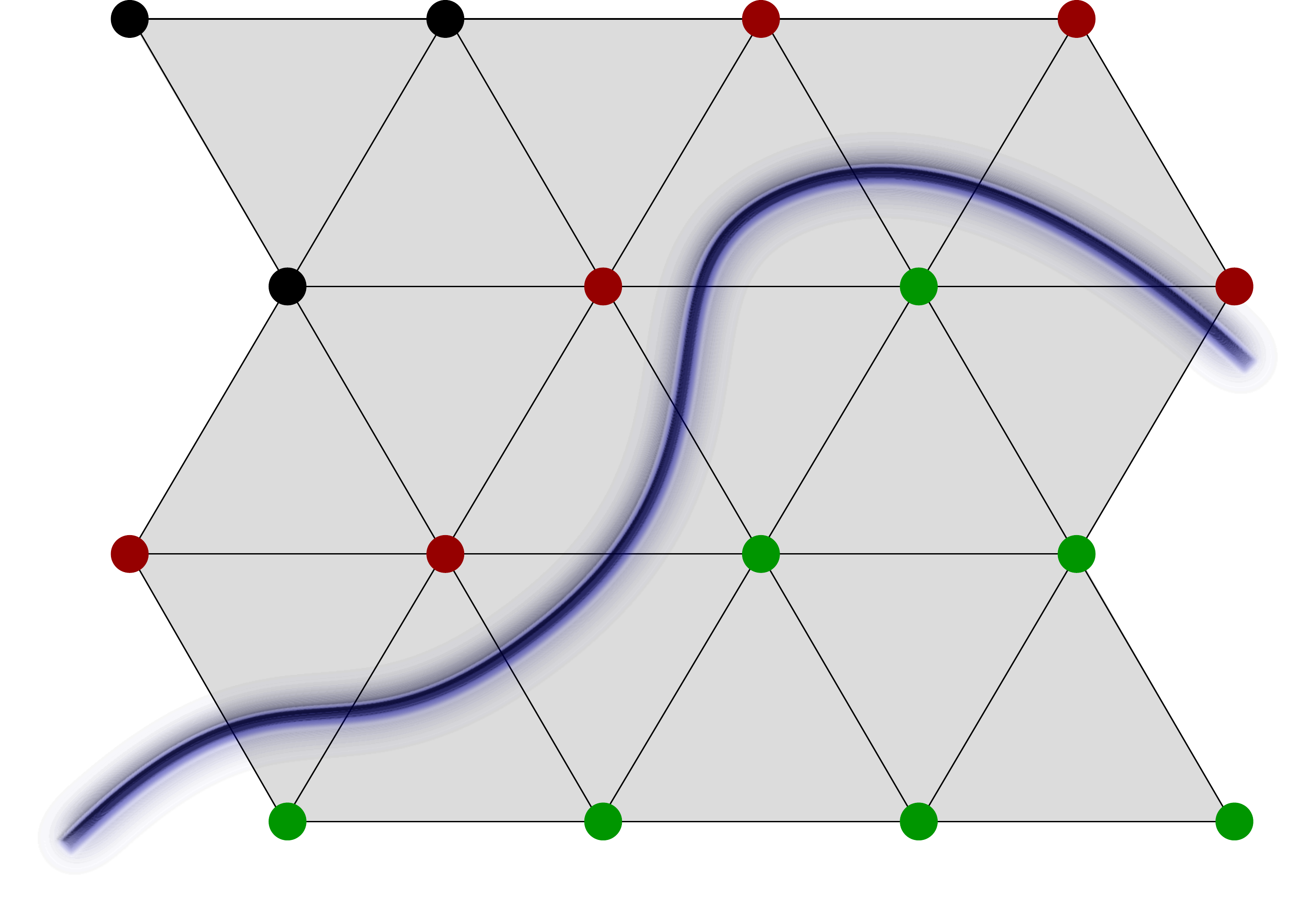}
  \end{center}
  \caption{Left: in green, nodes of $\mathcal{N}_\mathcal{P}$ and in black, nodes of $\mathcal{N}_\mathcal{Q}$. Right: In red are highligted nodes selected to belong to $\mathcal{N}(\Gamma^h_{n+1})$. The blurred blue line represent approximate location of $\Gamma^h_{n+1}$.}
  \label{fig:nodeSelec}
\end{figure}

For a given solution $(U_{n+1}, \vec X_{n+1})$, it is possible to split the set of nodes $\mathcal{N}$ into two subsets $\mathcal{N}_\mathcal{P}$ and $\mathcal{N}_\mathcal{Q}$:
\begin{equation}
  \left\{
  \begin{array}{lll}
    \mathcal{N}_\mathcal{P} & = & \{i\in\mathcal{N} \mid U_{n+1}^i > 0\}\\
    \mathcal{N}_\mathcal{Q} & = & \mathcal{N}\setminus\mathcal{N}_\mathcal{P}\\
  \end{array}
  \right.
\end{equation}
$\mathcal{N}_\mathcal{P}$ corresponds to the set of nodes which are inside \emph{the phase} at $t_{n+1}$ and $\mathcal{N}_\mathcal{Q}$ to the set of nodes which are outside \emph{the phase}. Interface $\Gamma^h_{n+1}$ is localized inside triangles of $\mathcal{T}_{n+1}$ having vertices belonging to $\mathcal{N}_\mathcal{P}$ and vertices belonging to $\mathcal{N}_\mathcal{Q}$, but the exact location is not known a priori. We now define:
\begin{equation}
  \mathcal{N}(\Gamma^h_{n+1}) = \{ i\in\mathcal{N}_\mathcal{Q} \mid i \text{ is connected to at least one node of } \mathcal{N}_\mathcal{P}\}
\end{equation}
$\mathcal{N}(\Gamma^h_{n+1})$ is the set of nodes which will be used to describe $\Gamma^h_{n+1}$ (Figure~\ref{fig:nodeSelec}). It is defined only once for time step $t_{n+1}$ and will remain unchanged during the successive updates of $\mathcal{T}_{n+1}$ nodes location. This step is called $localize\_interface$ in Algorithm \ref{alg:xm}. During the mesh update procedure, nodes from $\mathcal{N}(\Gamma^h_{n+1})$ will be moved to $\Gamma^h_{n+1}$ estimated location.

It is important to note that moving all nodes from $\mathcal{N}(\Gamma^h_{n+1})$ to $\Gamma^h_{n+1}$ is an ill-posed non-linear problem. The approach adopted is to compute and update an estimated optimal position for each node of $\mathcal{N}(\Gamma^h_{n+1})$ independently from each other and reiterate this process until convergence.

Let's consider $p\in\mathcal{N}(\Gamma^h_{n+1})$ a node which has been flagged to be belonging to the interface $\Gamma^h_{n+1}$. In order to obtain a good estimation of the interface $\Gamma^h_{n+1}$ location, we are considering the current solution $U_{n+1}$ only on the patch of triangles $\mathcal{T}^p_{n+1}$ sharing $p$ as a vertex. We denote $\mathcal{N}^p$ the set of nodes of $\mathcal{T}^p_{n+1}$ and $\bar{\mathcal{N}}^p = \mathcal{N}^p \setminus \mathcal{N}(\Gamma^h_{n+1})$ the set of nodes connected to $p$ which does not belong to $\mathcal{N}(\Gamma^h_{n+1})$. Finally, we define $\tilde{U}_{n+1}$ such as:
\begin{equation}
  \left\{
  \begin{array}{l}
    \tilde{U}_{n+1}^i = U_{n+1}^i\text{, }\forall i\in\mathcal{N}^p\text{, }i\neq p\\
    \tilde{U}_{n+1}^p = 0\\
  \end{array}
  \right.
  \label{eq:patchSolVal}
\end{equation}
Considering definition Eq.~\eqref{eq:defr}, $p$ belonging to the interface $\Gamma^h_{n+1}$ is equivalent to:
\begin{equation}
  r(\phi_p, \tilde{U}_{n+1},\vec{X}_{n+1}, U_{n},\vec{X}_{n}) = 0
  \label{eq:nodeBelongsInterface}
\end{equation}

We define now  a vector $\vec{v}\in\mathbb{R}^d$ and a new triangulation $\tilde{\mathcal{T}}_{n+1}^\vec{v}$ obtained from the triangulation $\mathcal{T}_{n+1}$ where the node $p$ is translated by $\vec{v}$. The problem of finding a new location for $p$ such as $p\in\Gamma^h_{n+1}$ can then be formulated as:
\begin{equation}
    \text{Find } \vec{v}\in\mathbb{R}^d \text{ such as } r(\phi_p, \tilde{U}_{n+1},\tilde{\vec{X}}_{n+1}^\vec{v}, U_{n},\vec{X}_{n}) = 0
  \label{eq:pbLocInterface}
\end{equation}
This problem does not have a unique solution. In order to simplify the search for a valid $\vec{v}$, we are only considering translations along edges connecting $p$ to $\bar{\mathcal{N}}^p$. Therefore, we only consider vectors $\vec{v} = \eta_{pj}\vec{v}_{pj}$ where:
\begin{equation}
  \left\{
  \begin{array}{l}
    \vec{v}_{pj} = \vec{X}_{n+1}(j) - \vec{X}_{n+1}(p)\text{, }j\in\bar{\mathcal{N}}^p\\
    \eta_{pj}\in[0, 1]\\
  \end{array}
  \right.
  \label{eq:defVsearch}
\end{equation}
In doing so, solving problem Eq.~\eqref{eq:pbLocInterface} is reduced to solving $|\bar{\mathcal{N}}^p|$ sub-problems:
\begin{equation}
  \left\{
  \begin{array}{l}
    \text{With $p$, $\tilde{U}_{n+1}$, $U_{n}$, $\vec{X}_{n}$, $\vec{v}_{pj}$ fixed, and $j\in\bar{\mathcal{N}}^p$,}\\
    \text{find } \eta_{pj}\in[0,1] \text{ such as } \\
    r(\eta_{pj}) =  r(\phi_p, \tilde{U}_{n+1},\tilde{\vec{X}}_{n+1}^{\eta_{pj}\vec{v}_{pj}}, U_{n},\vec{X}_{n}) = 0\\
  \end{array}
  \right.
  \label{eq:pbLocInterfaceSub}
\end{equation}

\begin{figure}[h!]
  \begin{center}
    \includegraphics[width=0.48\textwidth]{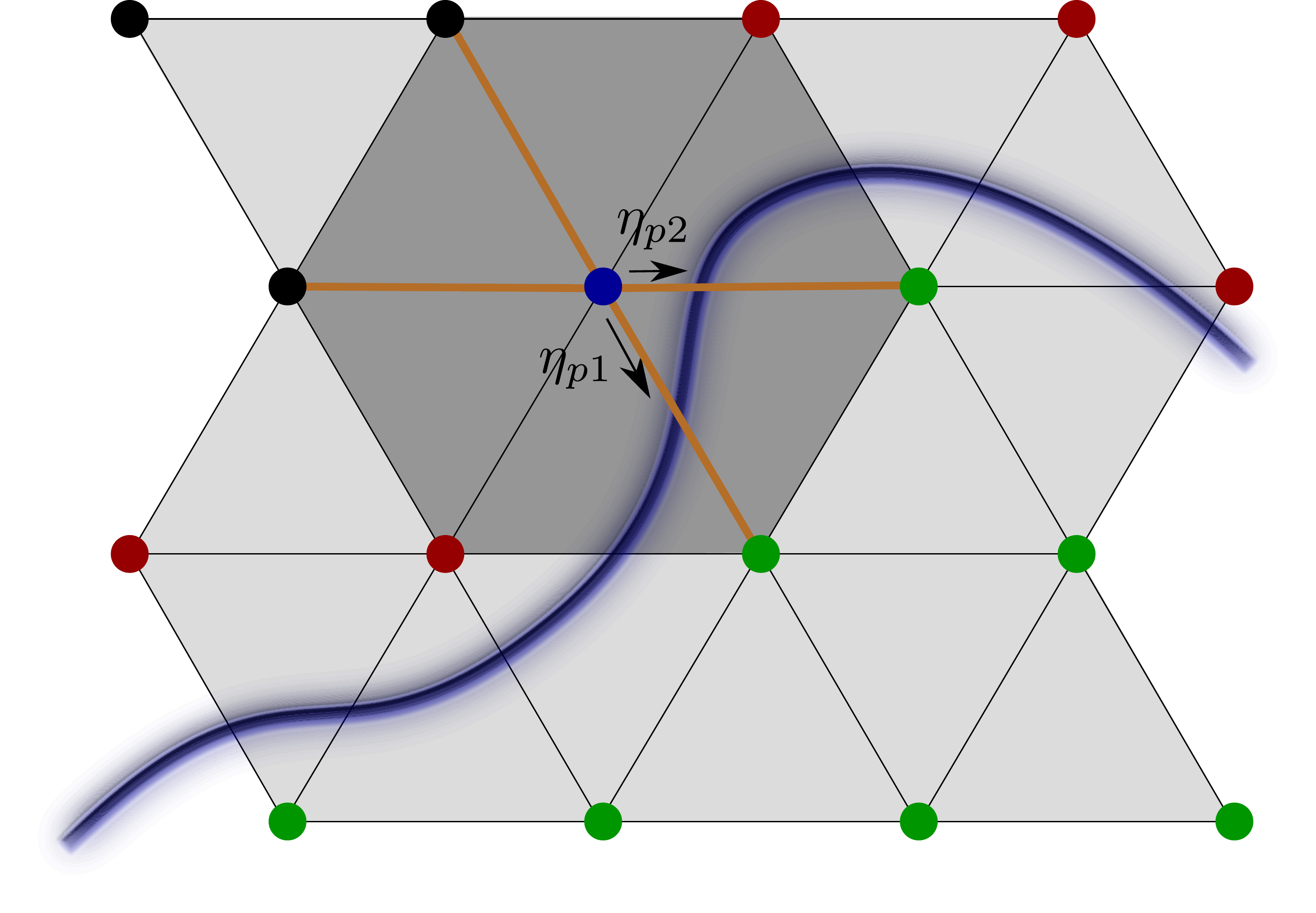}
    \hfill
    \includegraphics[width=0.48\textwidth]{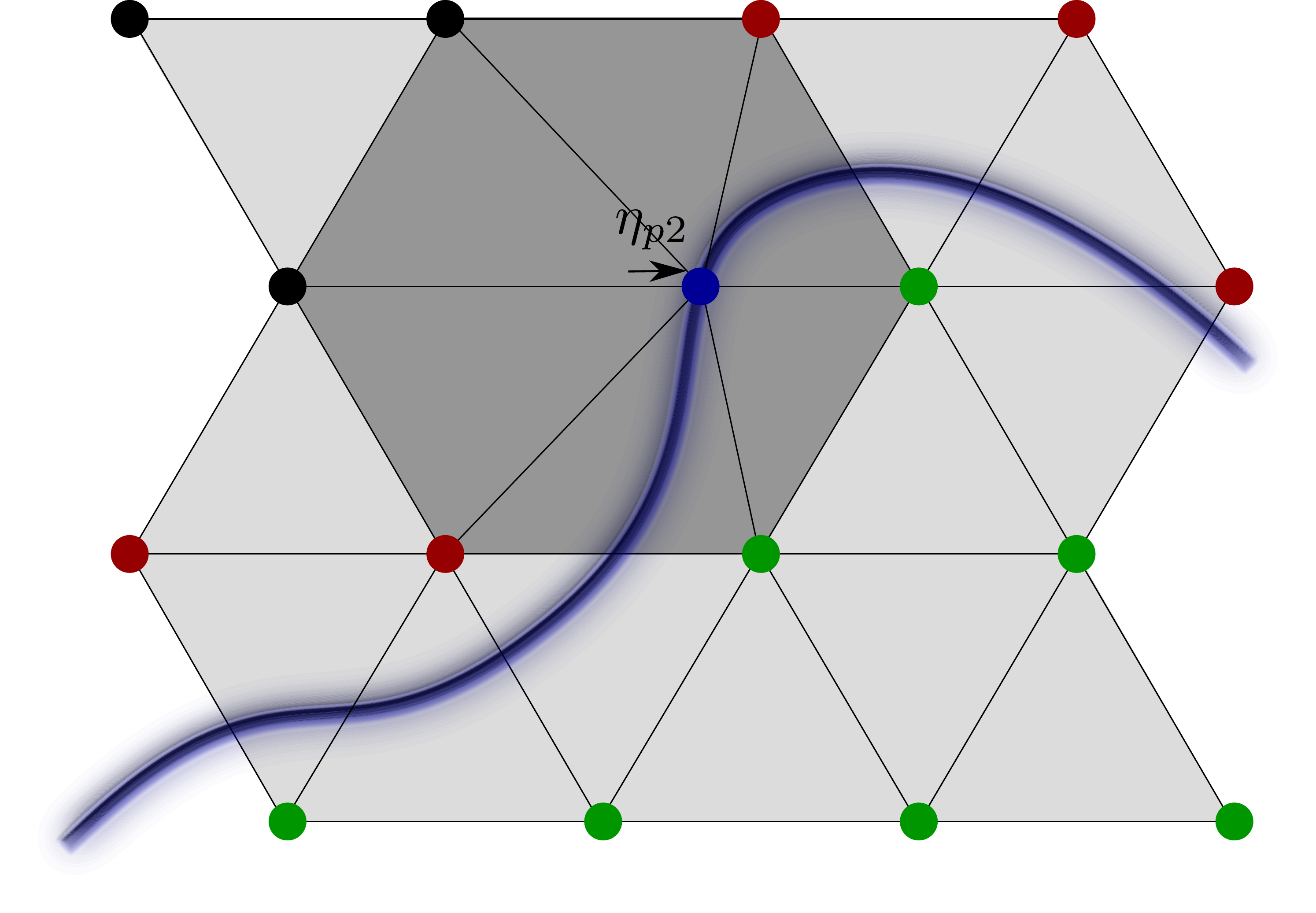}
  \end{center}
  \caption{Left: The node $p$ flagged to be moved onto estimated interface $\Gamma^h_{n+1}$ location is represented in blue. Dark grey area corresponds to the associated triangles patch $\mathcal{T}^p_{n+1}$. Orange line segments represent edges along which a translation of node $p$ is considered. $\eta_{p1}$ and $\eta_{p2}$ are the two displacement obtained verifying problem Eq.~\eqref{eq:pbLocInterfaceSub}. Right: Mesh obtained after choosing the displacement inducing the smallest mesh deformation.}
  \label{fig:nodeMovement}
\end{figure}

The sub-problem is not guaranteed to have a solution, yet given the strict monotonicity of $r(\eta_{pj})$ with respect to $\eta_{pj}$, any existing solution would be unique and determined using a bisection method.
Defining $\mathcal{S}_p = \{ (\eta_{pj}, \vec{v}_{pj}) \mid (\eta_{pj}, \vec{v}_{pj})\text{ verifies \ref{eq:pbLocInterfaceSub}} \}$ the set of possible translations to move node $p$ onto the estimated $\Gamma^h_{n+1}$, we choose the one inducing the smallest mesh deformation possible, that is to say the one for which $\|\eta_{pj}\vec{v}_{pj}\|$ is minimal. The corresponding node $j$ for which $\|\eta_{pj}\vec{v}_{pj}\|$ is minimal will be called node $p$ target and is noted $q_p$:
\begin{equation}
  q_p = \argmin_{\{j\mid\eta_{pj}\vec{v}_{pj}\in\mathcal{S}_p\}} \|\eta_{pj}\vec{v}_{pj}\|
  \label{eq:qp}
\end{equation}
The method described is illustrated in Figure~\ref{fig:nodeMovement}.

As mentioned previously $\Gamma^h_{n+1}$ actually depends on $U_{n+1}$, $\vec{X}_{n+1}$, $U_{n}$ and $\vec{X}_{n}$. Moving a node $p$ of $\mathcal{T}_{n+1}$ will modify the location of $\Gamma^h_{n+1}$ in its vicinity. Estimate obained for a point $p$ from solving Eq.~\eqref{eq:pbLocInterfaceSub} are no longer valid if $p$ neighbours are moved too. In practice, the translation $\vec v = \eta_{pq_p}\vec{v}_{pq_p}$ obtained is capped by a chosen value $\alpha$ to ease convergence:  $\vec v = \min(\eta_{pq_p}, \alpha)\vec{v}_{pq_p}$. We found $\alpha = 0.6$ to be working well.

\begin{figure}[h!]
  \begin{center}
    \includegraphics[width=0.75\textwidth]{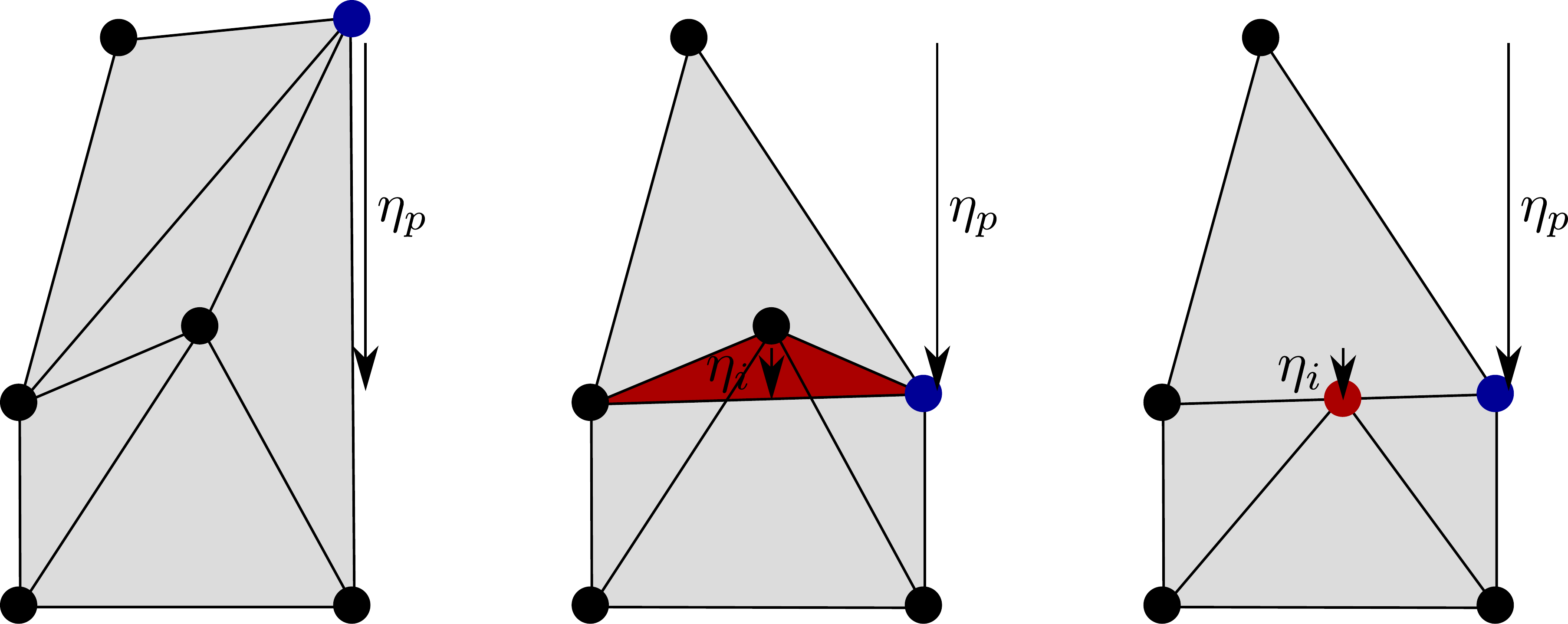}
  \end{center}
  \caption{Example of inverted element creation and repairing which can occur during the mesh update procedure. Left: mesh to be updated applying $\eta_p$ translation to the node highlighted in blue. Middle: this translation generates an inverted triangle depicted in red. This inverted triangle is flattened by applying translation $\eta_i$ to the node depicted in red. Right: resulting mesh after the $fix\_inverted\_elements$ procedure.}
  \label{fig:invertElem}
\end{figure}

Ultimately, the translation $\vec{v} = \min(\eta_{pq_p}, \alpha)\vec{v}_{pq_p}$ may result in the creation of inverted elements. If such a situation arises, node $p$ will be translated, and any resulting inverted elements will then undergo a corrective procedure referred to as $fix\_inverted\_elements$ in Algorithm \ref{alg:updateInterface} and illustrated in Figure~\ref{fig:invertElem}. This procedure aims to rectify the inverted elements by flattening them, producing elements with null or quasi-null measure. These elements are typically deemed unsuitable for finite element simulations due to their non-compliance with the convergence conditions outlined in \cite{babuvska1976angle}. However, as detailed in \cite{hannukainen2012maximum, duprez2019finite}, convergence of finite element solutions can still be achieved without strictly adhering to these conditions. A more relaxed criterion is developed in \cite{kuvcera2016necessary}. Additionally, as outlined in \cite{duprez2019finite}, acceptable patterns of these poorly shaped elements only affect the conditioning of finite element matrices. In the work presented here, all matrix systems are solved with direct solvers, effectively eliminating any concerns regarding the poor conditioning of operators. In cases where an iterative solver for matrix systems is required (e.g., for performance purposes with 3D problems), we recommend the approach proposed in \cite{moes2023extreme}, which involves bounding from below the determinant of the Jacobian of very small elements.

\begin{algorithm}[h!]
  \caption{$update\_interface$ method.}
  \textbf{Input:} $\mathcal{N}(\Gamma^h_{n+1}), (U_n, \vec{X}_n) , (U_{n+1}^{k}, \vec{X}_{n+1}^{k})$\;
  $\alpha = 0.6$\;
  \For{$p\in\mathcal{N}(\Gamma^h_{n+1})$}{
    Find $\mathcal{S}_p$ solving problem \ref{eq:pbLocInterfaceSub}\;
    \If{$\mathcal{S}_p = \emptyset$}{
      $\vec{v}[p] = \vec{0}$\;
    }
    \Else{
      $\vec{v}[p] = \displaystyle\argmin_{\vec{u}\in\mathcal{S}_p} \|\vec{u}\|$\;
    }
  }
  \For{$p\in\mathcal{N}$}{
    $\vec{X}^{k+1}_{n+1}(p) = \vec{X}^{k}_{n+1}(p)$\;
  }
  \For{$p\in\mathcal{N}(\Gamma^h_{n+1})$}{
    $\vec{X}^{k+1}_{n+1}(p) = \vec{X}^{k+1}_{n+1}(p) + \vec{v}[p]$\;
    $\vec{X}_{n+1}^{k+1} \gets fix\_inverted\_elements(\vec{X}_{n+1}^{k+1})$\;
  }
  \textbf{Output:} $\vec{X}^{k+1}_{n+1}$\;
  \label{alg:updateInterface}
\end{algorithm}

\subsection{Mesh velocity definition}
\label{sec:meshVel}
The definition of the mesh velocity $\vec{w}$ within a time interval $[t_n, t_{n+1}]$ remains to be addressed. We know that $\vec{w}$ satisfies:
\begin{equation}
  \dint_{t_n}^{t_{n+1}} \vec{w}(\vec{X}_0, t) \,\text{d}t = \vec{X}(\vec{X}_0, t_{n+1}) - \vec{X}(\vec{X}_0, t_{n}) = \vec{X}_{n+1} (\vec{X}_0) - \vec{X}_{n} (\vec{X}_0)\text{, }\forall \vec{X}_0\in\Omega
  \label{eq:constrW}
\end{equation}
We have the flexibility to choose any $\vec{w}$ satisfying Eq.~\eqref{eq:constrW}. The natural modeling choice proposed in \cite{moes2023extreme} is a mesh velocity $\vec{w}$ such as $\left.\dpfrac{\vec w}{t}\right|_{\vec X} = \vec 0$ in $[t_n, t_{n+1}]$. Using definition established in Eq.~\eqref{eq:defXp}, it writes:
\begin{equation}
  \left\{
  \begin{array}{l}
    \vspace{0.5em}
    \vec{w}(\vec X) = \dsum_{i\in\mathcal{N}} \vec W^i \phi_i(\vec{X})\\
    \vec{W}^i = \dfrac{\vec{X}_{n+1}(i) - \vec{X}_{n}(i)}{t_{n+1}-t_n}\\
  \end{array}
  \right.
  \label{eq:defaultW}
\end{equation}

\begin{figure}[h!]
  \begin{center}
    \includegraphics[width=0.43\textwidth]{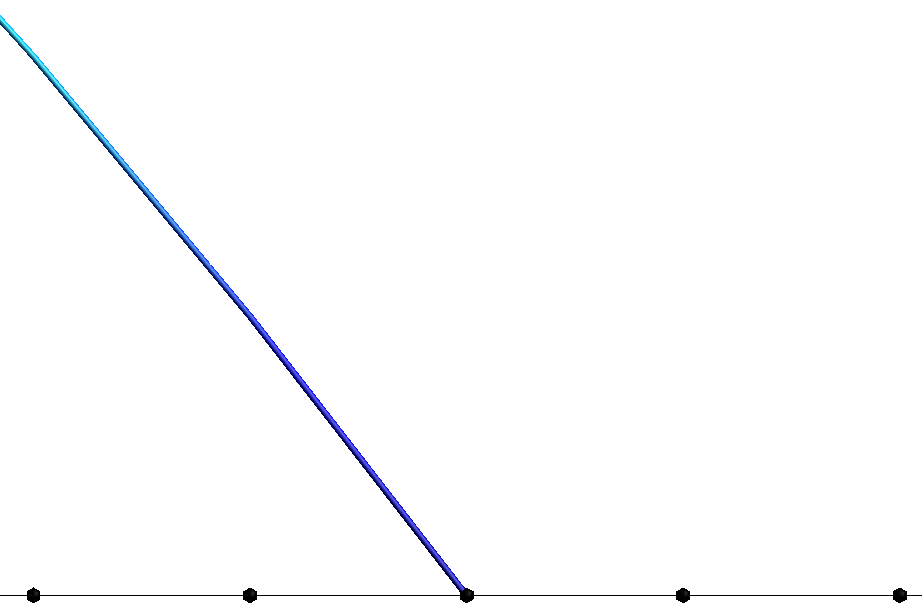}
    \hspace{0.03\textwidth}
    \includegraphics[width=0.43\textwidth]{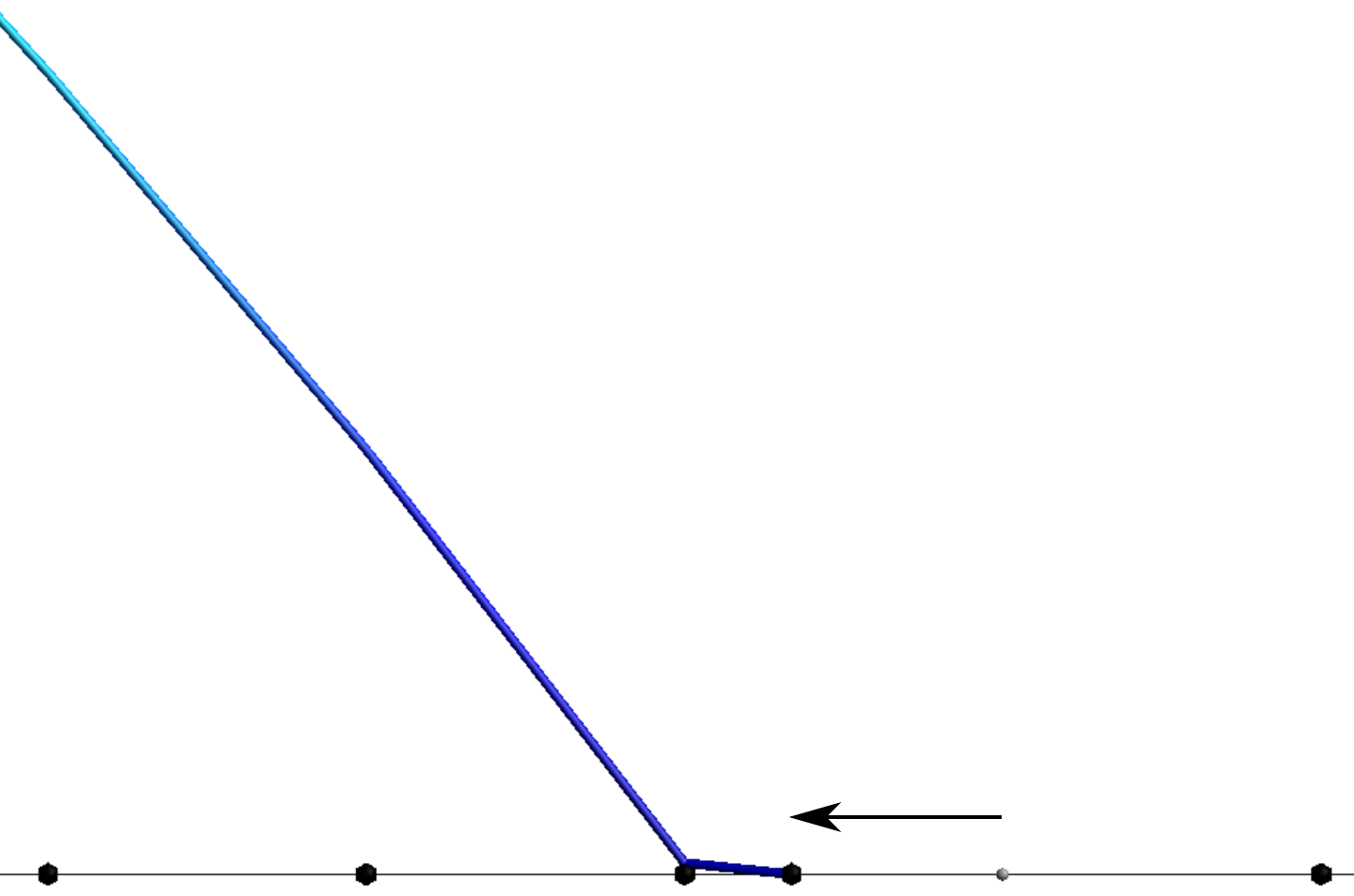}
  \end{center}
  \caption{Colored line represent the solution $u$ in the phase $\mathcal{P}_u$ and black spheres correspond to mesh nodes. Left: mesh and solution at $t_n$. Right: converged mesh and solution obtained at $t_{n+1}$ when $\vec w$ is defined using Eq.~\eqref{eq:defaultW}. Arrow represents the interface node's displacement.}
  \label{fig:oldW}
\end{figure}

This option is suitable for addressing \twop\, problems, however it is not appropriate for \onep\, problems due to its inadequate representation of the underlying physics. An example of unwanted behavior using this mesh velocity modelization is illustrated Figure~\ref{fig:oldW}. The converged solution obtained after one time step shows that the \emph{phase} $\mathcal{P}_u$ spreads more than expected and on a very thin sheet. As this result is obtained with a fully converged solution (there is no Lagrange multipliers left and $R((U_{n}, \vec{X}_{n}), (U_{n+1}, \vec{X}_{n+1})) = 0$), this behavior is not the consequence of the resolution method but originates from a modelization problem.

\begin{figure}[h!]
  \begin{center}
    \includegraphics[width=0.43\textwidth]{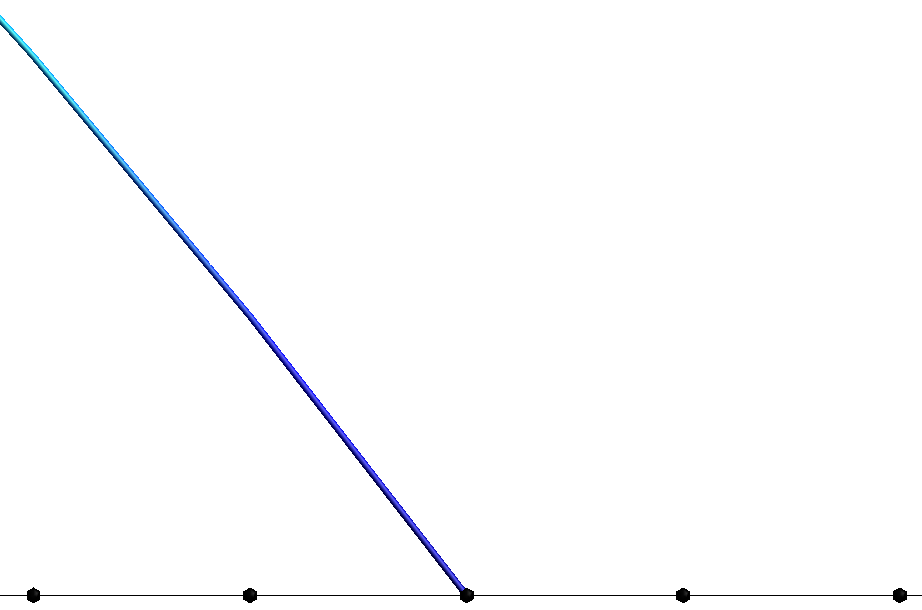}
    \hspace{0.03\textwidth}
    \includegraphics[width=0.43\textwidth]{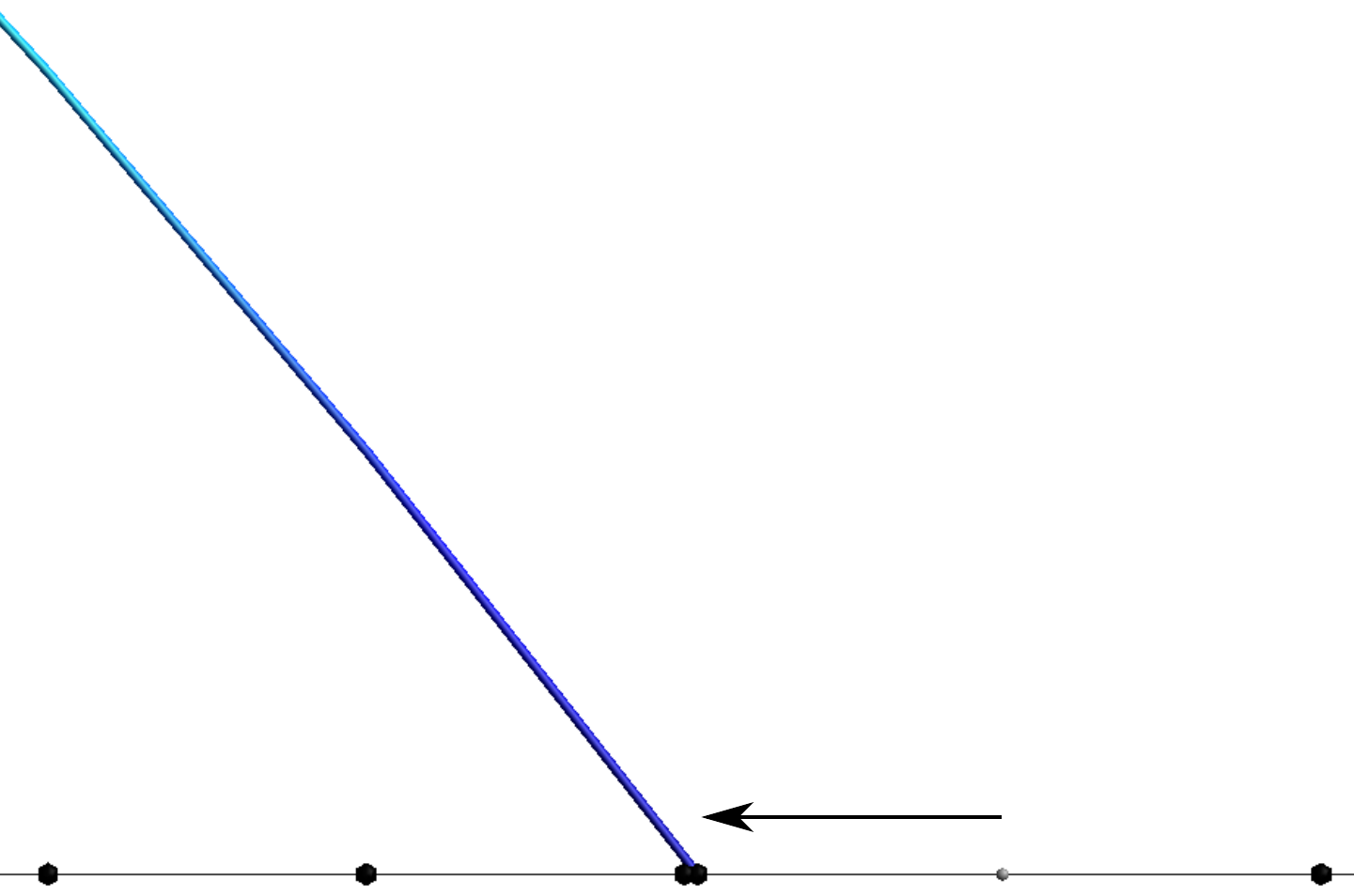}
  \end{center}
  \caption{Colored line represent the solution $u$ in the phase $\mathcal{P}_u$ and black spheres correspond to mesh nodes. Left: mesh and solution at $t_n$. Right: converged mesh and solution obtained at $t_{n+1}$ when $\vec w$ is defined using Eq.~\eqref{eq:newW}. Arrow represents the interface node's displacement.}
  \label{fig:newW}
\end{figure}

Therefore, the mesh velocity modelization choice made here is:
\begin{equation}
  \left\{
  \begin{array}{l}
    \vec{W}^p(t) = \dfrac{\vec{X}_{n+1}(p) - \vec{X}_{n}(q^0_p)}{t_{n+1}-t_n}\text{, if $p\in\mathcal{Q}_{u^h_n}$}\\
    \vec{W}^p(t) = \dfrac{\vec{X}_{n+1}(p) - \vec{X}_{n}(p)}{t_{n+1}-t_n}\text{, else.}\\
  \end{array}
  \right.
  \label{eq:newW}
\end{equation}
where $q^0_p$ is the target obtained for $p$ on the initial guess $\mathcal{T}^0_{n+1}$ solving Eq.~\eqref{eq:qp}. Further justifications of this choice can be found in \ref{ann:meshVel}. Results obtained using Eq.~\eqref{eq:newW} as the mesh velocity modeling are presented in Figure~\ref{fig:newW}. The unwanted spread is no longer present, and the behavior aligns more closely with expectations.

\section{Numerical results}
To evaluate the proposed numerical scheme's performance, we will perform a three-part assessment. First, we will compare the solution quality obtained with our scheme to that of the classical FEM for the Barenblatt-Pattle family of solutions. Second, we will investigate the scheme's ability to simulate waiting time phenomena. Finally, we will assess the method's robustness to changes in the phase topology.

\subsection{Barenblatt-Pattle solution}

This section leverages the Barenblatt-Pattle solutions defined in Eq.~\eqref{eq:barenblatt} to assess the accuracy of the numerical scheme. The domain $\Omega$ is discretized into a uniform mesh $\mathcal{T}_0$ with mesh size $h$. The numerical solution obtained on this discretization is denoted by $u^h$. Simulations are performed for different values of $m$ and $h$, and a backward Euler approximation in time ($\theta = 1$) is employed. The following parameter values are used:

\begin{equation}
\kappa = 1, \quad r = 0.15, \quad t_0 = (kr^2)^{\frac{1}{\alpha m + 2b}}, \quad C = (t_0)^{\alpha m}
\label{eq:valBPsim}
\end{equation}

To assess the quality of numerical solutions, we examine two metrics: the relative interface localization error $f_r(t)$ and the relative mass variation $M_r(t)$. The former, defined in Eq.~\eqref{eq:relativeIntLocError}, gauges interface position accuracy relative to the mesh size $h$ , where $f(t)$ represents the analytical interface position (Eq.~\eqref{eq:bpIntF}) and $\mathcal{N}(\Gamma(t))$ denotes nodes on the interface at time $t$. The latter, defined in Eq.~\eqref{eq:relativeMassVar}, evaluates mass preservation throughout the simulation. 

\begin{equation}
    f_r(t) = \frac{1}{|\mathcal{N}(\Gamma(t))|} \sum_{p \in \mathcal{N}(\Gamma(t))} \frac{|\| \vec{x}^p(t) \| - f(t)|}{h} \label{eq:relativeIntLocError}
\end{equation}

\begin{equation}
    M_r(t) = \frac{\int_\Omega u^h(t) \mathrm{d}\Omega - \int_\Omega u^h(t_0) \mathrm{d}\Omega}{\int_\Omega u^h(t_0) \mathrm{d}\Omega} \label{eq:relativeMassVar}
\end{equation}

\begin{figure}[h!]
\centering
  \includegraphics[width=0.48\linewidth]{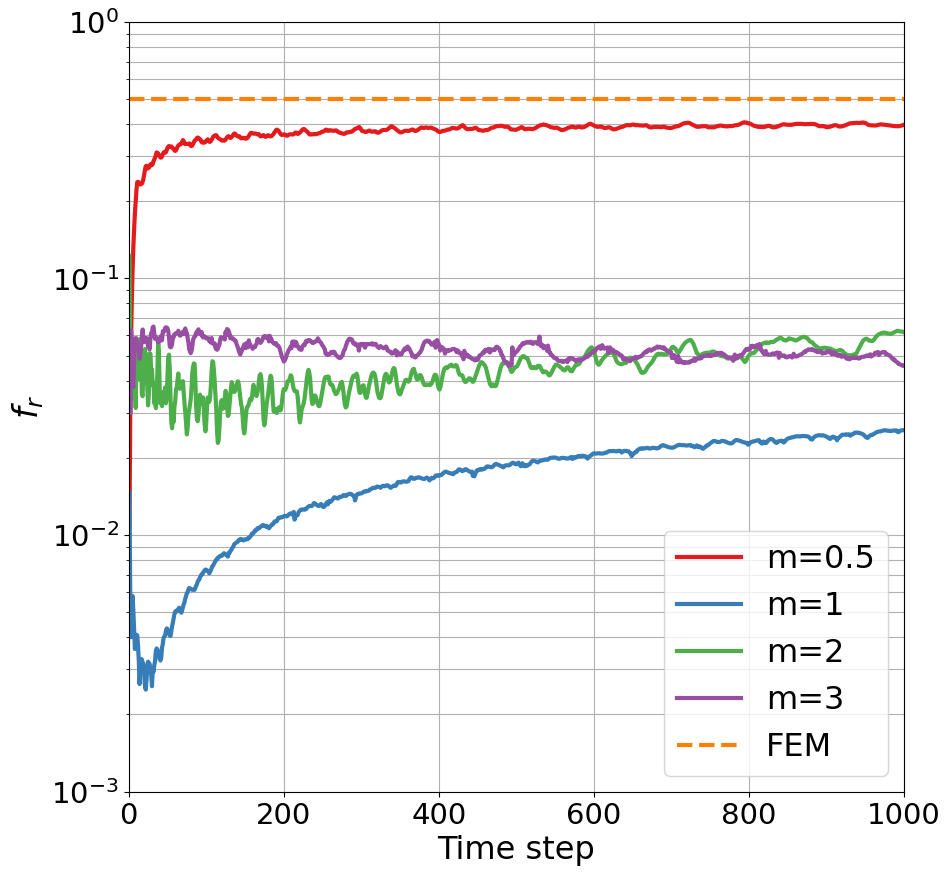}
  \includegraphics[width=0.48\textwidth]{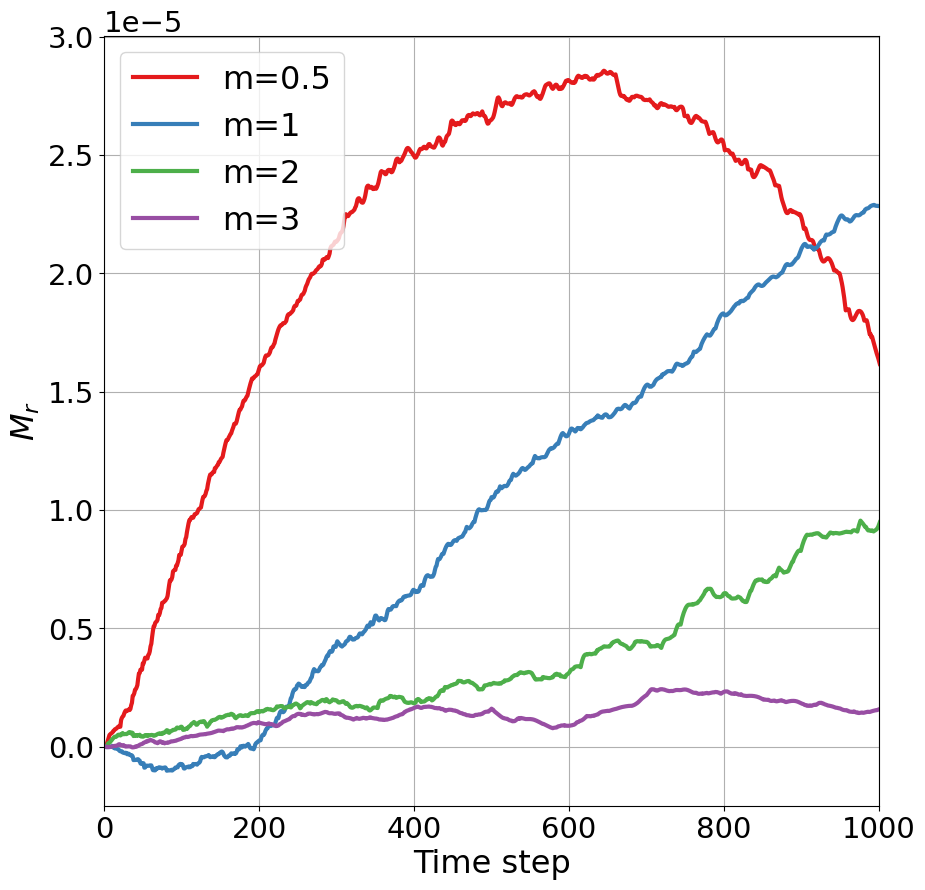}
\caption{Left: Relative interface localization error $f_r$ (Eq.~\eqref{eq:relativeIntLocError}) for various $m$ values with respect to the number of time iterations. The dotted orange line represents the average interface localization error relative to mesh size using the classical FEM. Right: Relative mass variation $M_r$ (Eq.~\eqref{eq:relativeMassVar}) with respect to the number of time iterations.}
\label{fig:barenblattMass}
\end{figure}

These two quantities are displayed on Figure~\ref{fig:barenblattMass} for various values of $m$. With an average relative interface localization error of $0.5$ for a classical FEM, the results show that the \xm{} approach achieves an order of magnitude higher precision in interface localization for $m \geq 1$. However, for $m<1$, no significant gain on interface localization is achieved. As it will be highligthed a bit further in this section, this is due to the fact that Barenblatt-Pattle solutions for $m<1$ are $\mathcal{C}^1(\Omega\times[t_0, T])$. The right plot on Figure~\ref{fig:barenblattMass} demonstrates that the mass conservation property, discussed in Section~\ref{sec:PME}, is maintained up to the user-specified \xm{} solver tolerance.

\begin{figure}[h!]
  \begin{center}
    \includegraphics[width=0.48\textwidth]{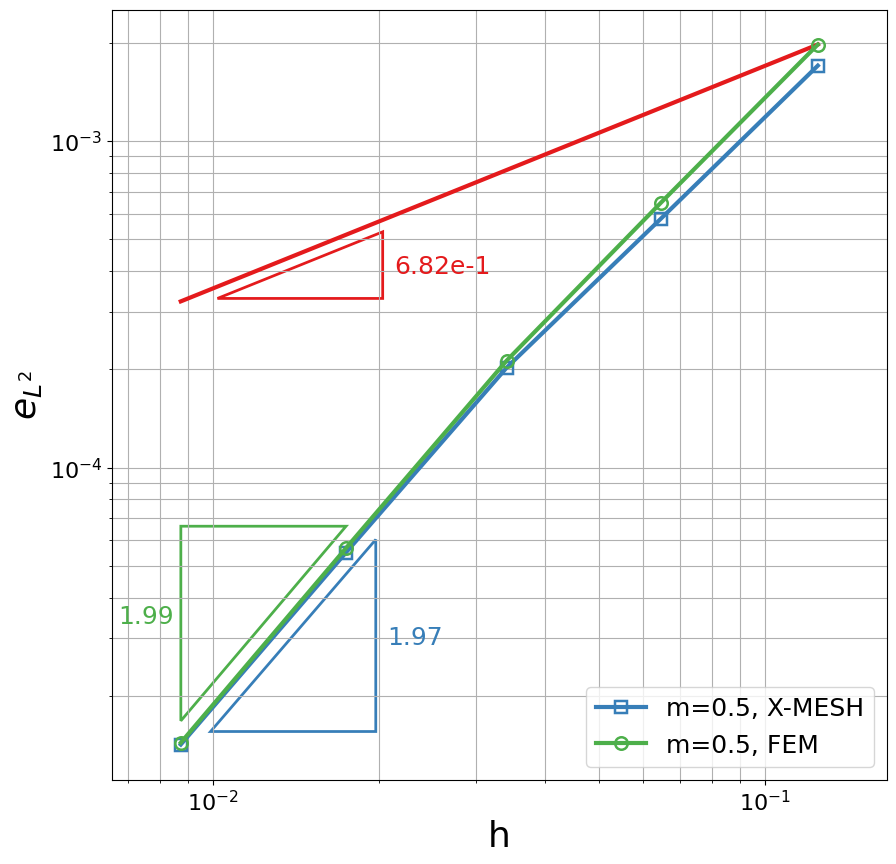}
    \includegraphics[width=0.48\textwidth]{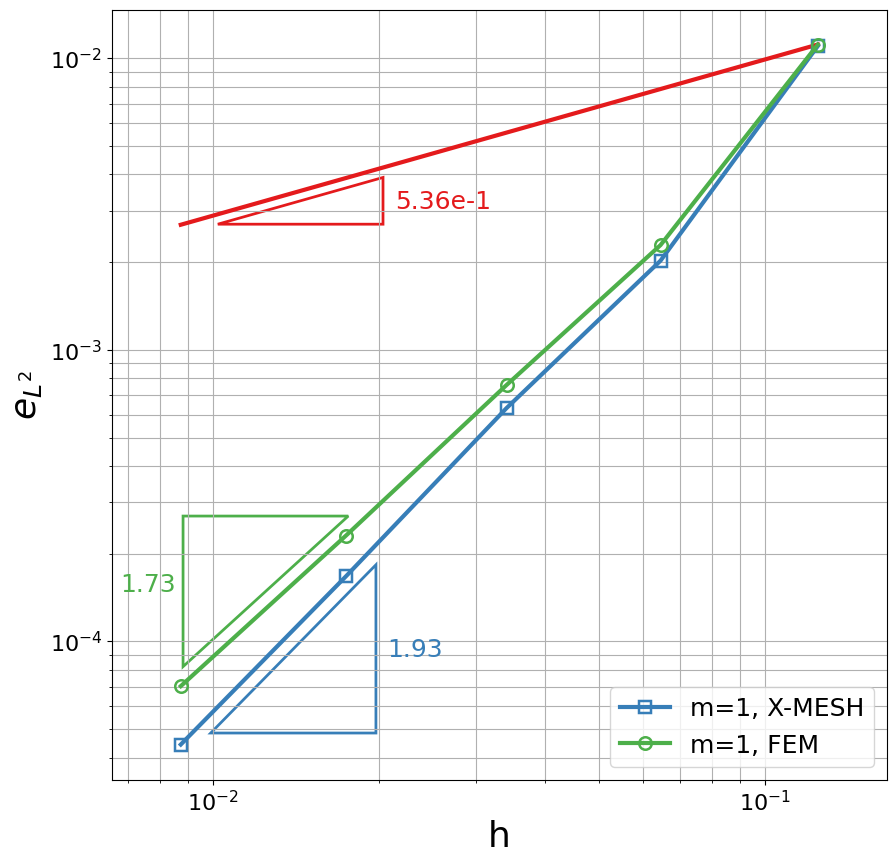}
    \includegraphics[width=0.48\textwidth]{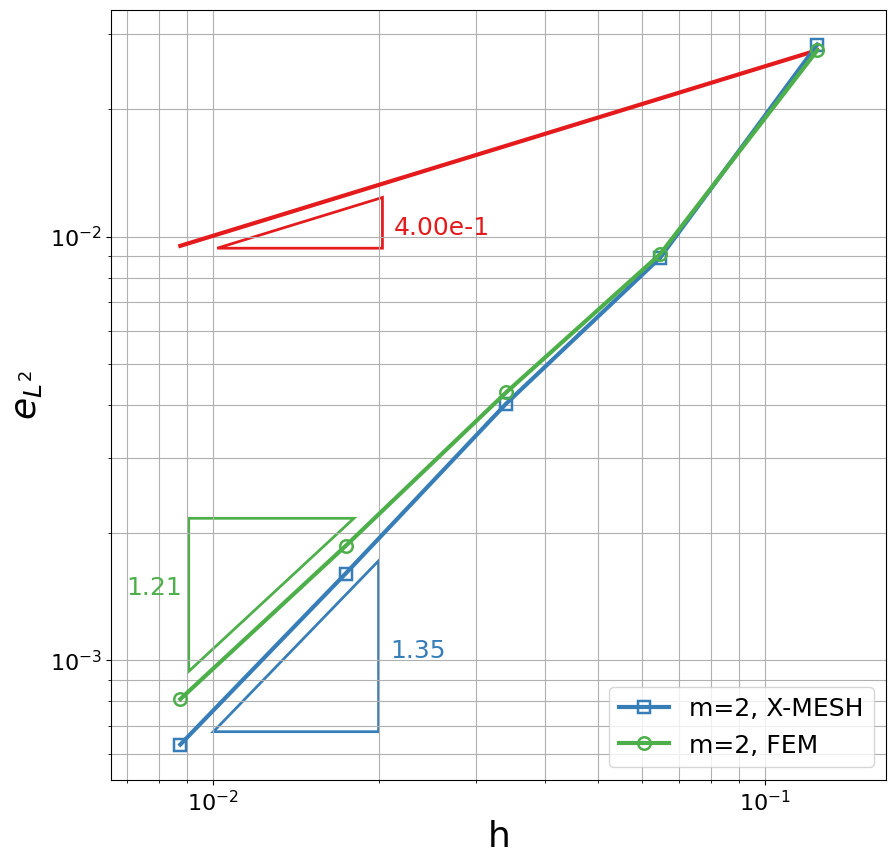}
    \includegraphics[width=0.48\textwidth]{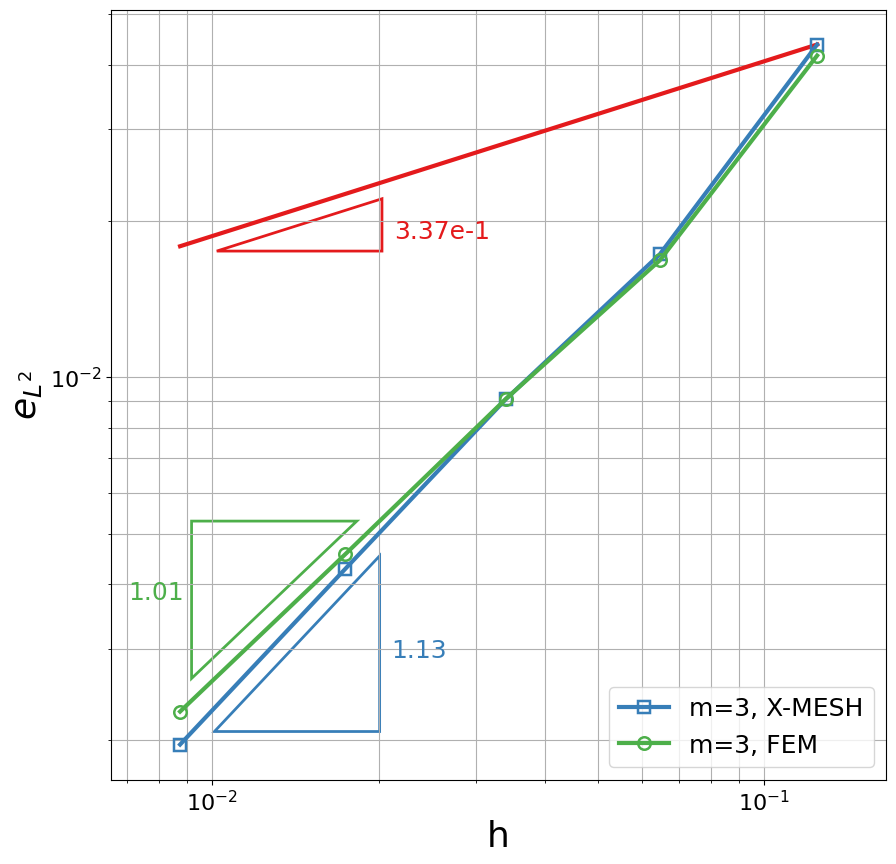}
  \end{center}
  \caption{$L^2$ error on Barenblatt-Paddle solutions for various values of $m$ regarding $h$. The solid red line correspond to the analytical $L^2$ error upper bound defined Eq.~\eqref{eq:convRate}.}
  \label{fig:barenblattPrec}
\end{figure}

Finally, we can analyse the behavior of the $L^2$ error between the numerical solution $u^h$ and the analytical solution $u_{BP}$ regarding $h$:
\begin{equation}
  e_{L^2}(u^h) = \|u^h - u_{BP}\|_{L^2(\Omega\times[t_0, T])} = \left(\dint_{t_0}^T \dint_\Omega (u^h-u_{BP})^2\id\Omega \id t\right)^{\frac{1}{2}},
  \label{eq:l2norm}
\end{equation}
The first thing to note from the results shown in Figure~\ref{fig:barenblattPrec} is that both the classical FEM and the \xm{} approach have a better convergence rate than the upper bound Eq.~\eqref{eq:convRate}. This is due to the overall smoothness of the Barenblatt-Pattle solutions.
Then, results can be divided in two categories: results for $m<1$ and $m\geq 1$. For the case $m<1$, it is known that solutions belong to $\mathcal{C}^1(\Omega\times [t_0, T])$. As a result, FEM methods will have a convergence rate in $\mathcal{O}(h^2)$ with linear basis functions and therefore cannot be expected to be better with the \xm{} approach as the only difference between both is the location of nodes belonging to the interface $\Gamma$. As we can see of the top left plot on Figure~\ref{fig:barenblattPrec}, \xm{} approach allows to obtain a slightly better solution regarding the $L^2$ norm than with a classical FEM method with coarse meshes, but converges asymptotically to the same error when $h\rightarrow 0$.

The situation is different for $m\geq 1$. Indeed, corresponding solutions are only $\mathcal{C}^0$ at the interface $\Gamma$. Figure~\ref{fig:barenblattPrec} shows that \xm{} approach leads to a convergence rate approximately $11\%$ greater than with a classical FEM method, highlighting that the interface $\Gamma$ being fully represented by the spatial discretization $\mathcal{T}$ at all time step allows to slightly reduce error commited during the numerical resolution. If a more significant convergence rate on the global $L^2$ error regarding the number of elements is sought-after, the \xm{} approach should be coupled with a method allowing to generate a finer discretization close to the interface $\Gamma$, as this is the region where strongest gradients can be found (for $m>1$, the solution's slope is infinite at $\Gamma$).

\begin{figure}[h!]
  \begin{center}
    \includegraphics[width=0.95\textwidth]{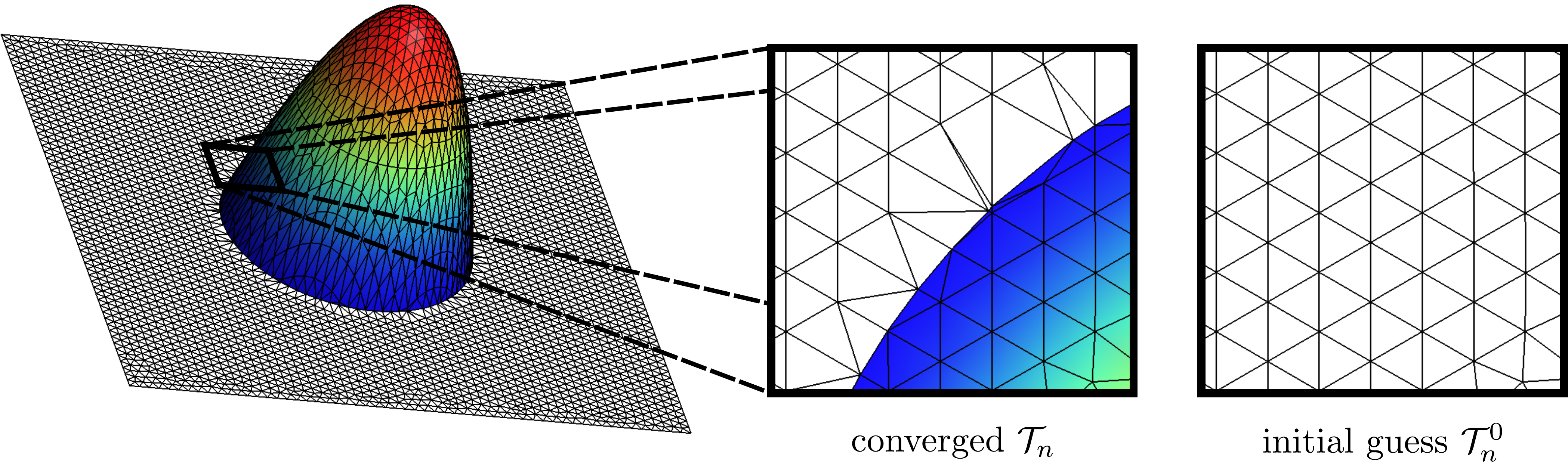}
  \end{center}
  \caption{Numerical solution for a given time step. On the right is the initial guess made for the mesh at $t_{n}$ and in the middle the converged mesh $\mathcal{T}_{n}$. Colored area corresponds to part of the domain where the solution is strictly positive.}
  \label{fig:exNumSol}
\end{figure}

Figure~\ref{fig:exNumSol} highlights the numerical solution's behavior at interface $\Gamma$. Mesh is adapted to fit interface $\Gamma$ method ensuring non-negativity of the solution.

\newpage
\subsection{Waiting time phenomena}

This section explores the behavior of the \xm{} approach when encountering initial conditions that result in a waiting time phenomena, as outlined in Section~\ref{sec:PME}, where the phenomena entails the interface $\Gamma(u(t))$ remaining stationary for a finite duration.

We specifically investigate an example detailed in \cite{aronson1983initially} for $d=1$ and $\kappa = 1$. The initial condition is:
\begin{equation}
  u_0(\vec{x}) =
  \left\{
  \begin{array}{ll}
    \left((1-\gamma)\sin^2(\vec{x}) + \gamma\sin^4(\vec{x})\right)^\frac{1}{m} & \text{, for }\vec{x}\in[-\pi, 0]\\
    0 & \text{, for }\vec{x}\in\mathbb{R}\setminus[-\pi, 0]\\
  \end{array}
  \right.
  \label{eq:waitingTime1}
\end{equation}
where an explicit expression for the waiting time $t_w$ is applicable when $\gamma\in[0, \frac{1}{4}]$, given by:
\begin{equation}
  t_w(m, \gamma) = \frac{1}{2(m+2)(1-\gamma)}
  \label{eq:wTanalytic}
\end{equation}
This analytical value serves as a benchmark to validate our numerical scheme. Simulations are performed across various values of $m$ and $\gamma = 0.05$.

\begin{figure}[h!]
  \begin{center}
    \includegraphics[width=0.48\textwidth]{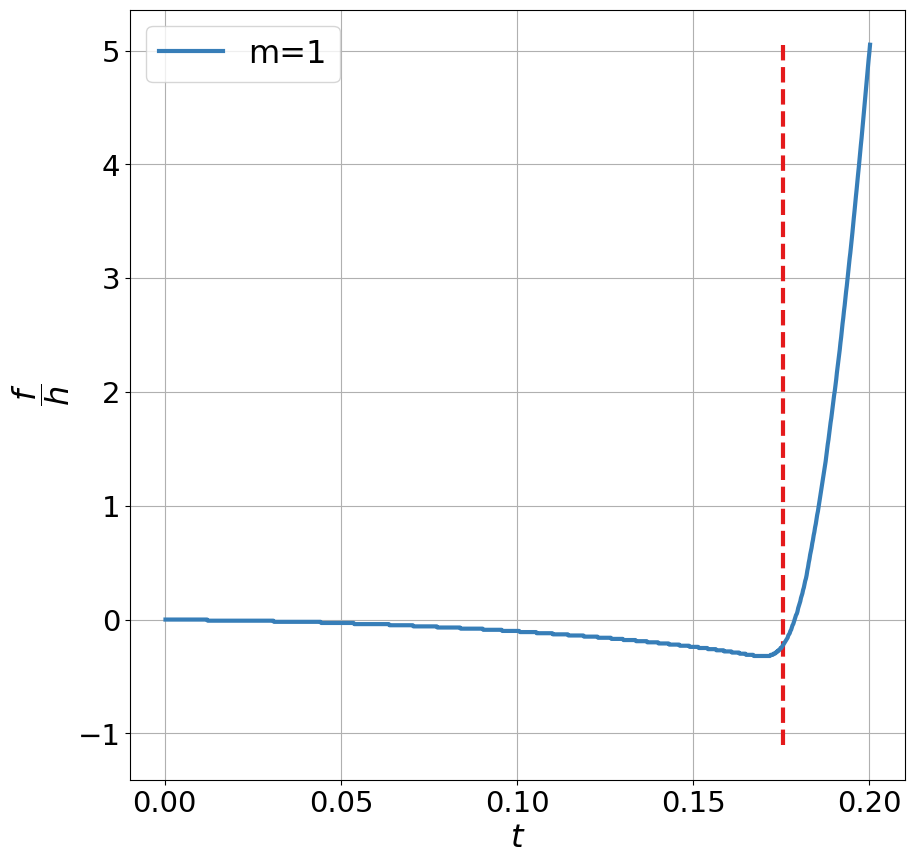}
    \includegraphics[width=0.48\textwidth]{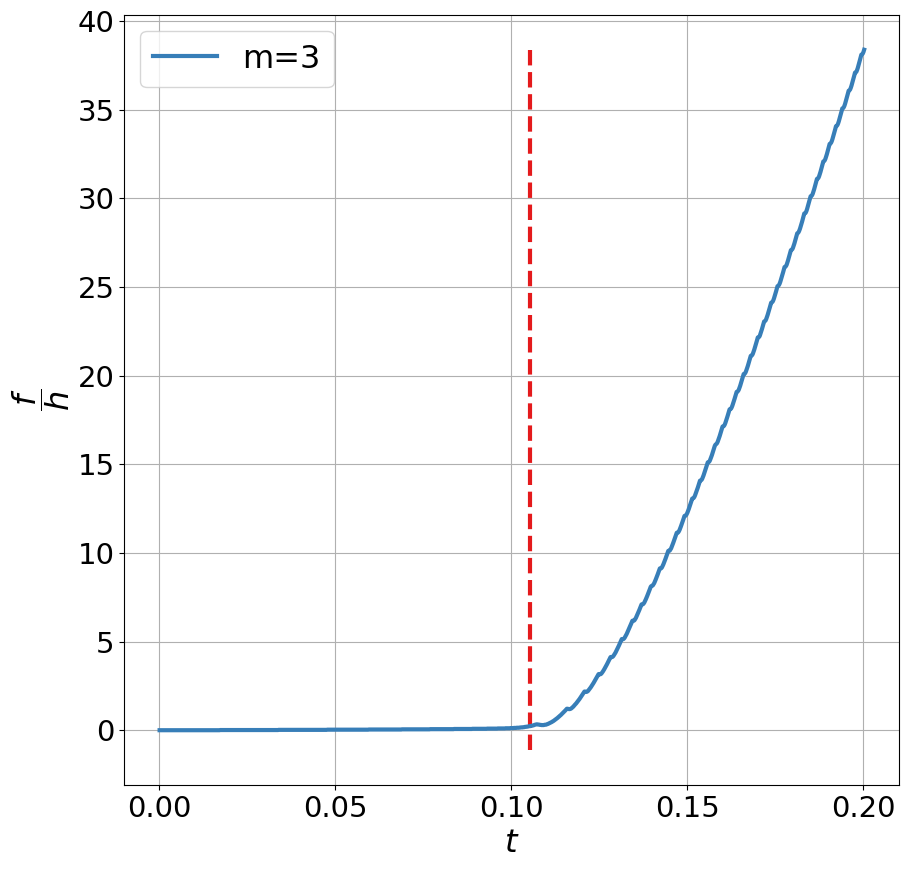}
    \includegraphics[width=0.48\textwidth]{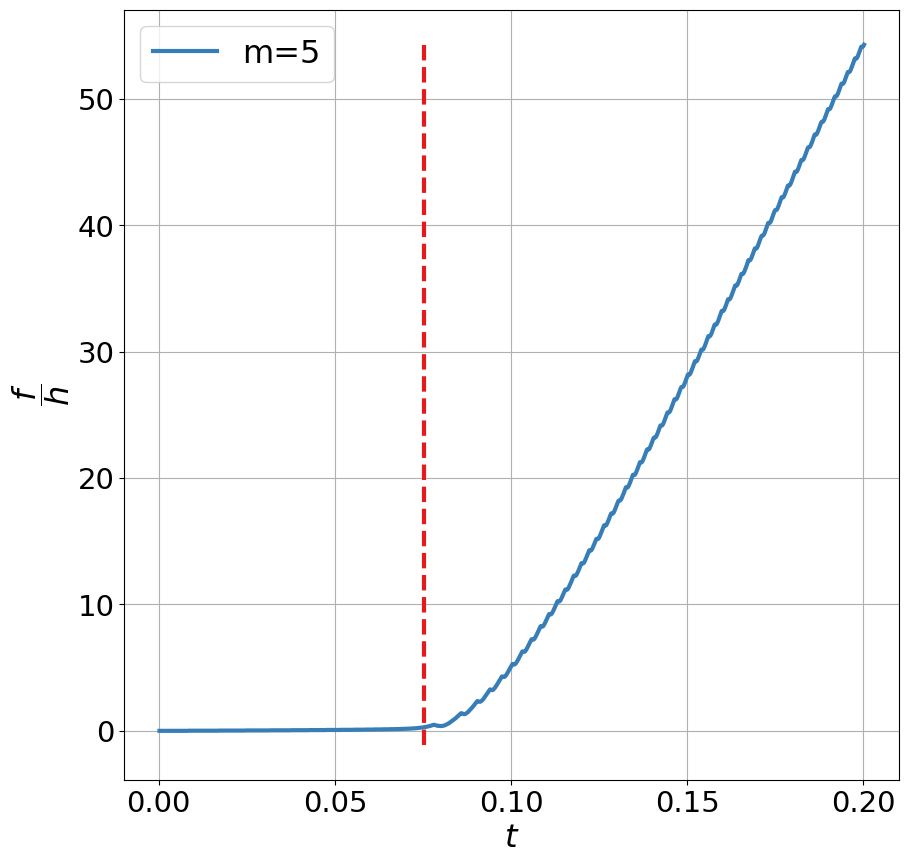}
    \includegraphics[width=0.48\textwidth]{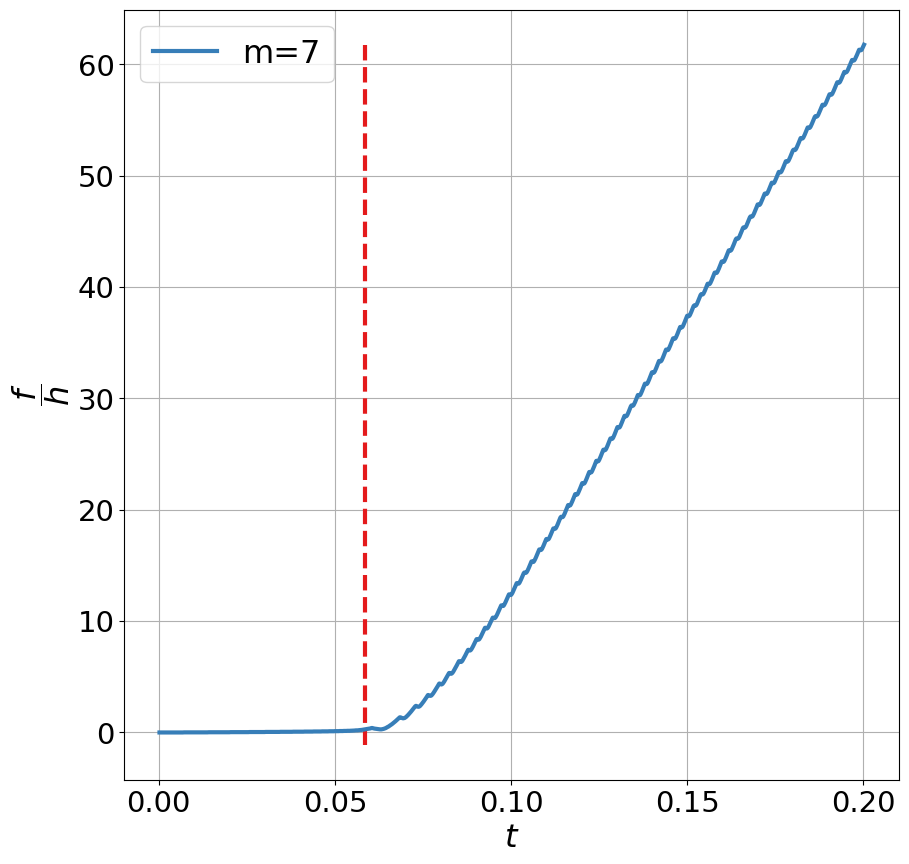}
  \end{center}
  \caption{Results obtained with the \xm{} approach for the initial condition Eq.~\eqref{eq:waitingTime1}. Continuous blue line: interface $\Gamma$ displacement $f$ relative to the mesh size $h$ regarding time. Dotted red line: analytical waiting time $t_w$ from \cite{aronson1983initially}.}
  \label{fig:wt1D}
\end{figure}

Figure~\ref{fig:wt1D} illustrates the movement of the interface $\Gamma$ over time $f(t)$ relative to the mesh size $h$. We observe that for all values of $m$, $\Gamma$ begins significant motion at $t_w(m)$. However, within the interval $[0, t_w(m)]$, $\Gamma$ exhibits some retrograde motion. While the magnitude of this motion is significantly lower than the motion observed for $t > t_w$, it contradicts the properties of the interface velocity of the PME's solutions (Eq.~\eqref{eq:intVelocity}). This discrepancy is particularly pronounced in the results for $m = 1$.

The waiting time phenomenon arises because $\Gamma$ remains stationary until the solution's slope near the interface becomes sufficiently steep to induce movement. On a fixed mesh, this can only be achieved through solution evolution at nodes, forcing a forward motion of the interface. In contrast, the \xm{} approach allows for this to be accomplished by moving nodes belonging to $\mathcal{N}(\Gamma)$. However, this flexibility introduces the possibility of capturing the unphysical retrograde motion observed here. Importantly, the magnitude of this unphysical motion is capped by the mesh size $h$.

Putting aside the capped, unphysical retrograde movement, waiting times obtained are in good agreement with analytical solutions (Eq.~\eqref{eq:wTanalytic}).

\newpage
\subsection{Interface Topology Changes}

\begin{figure}[h!]
  \begin{center}
    \includegraphics[width=0.24\textwidth]{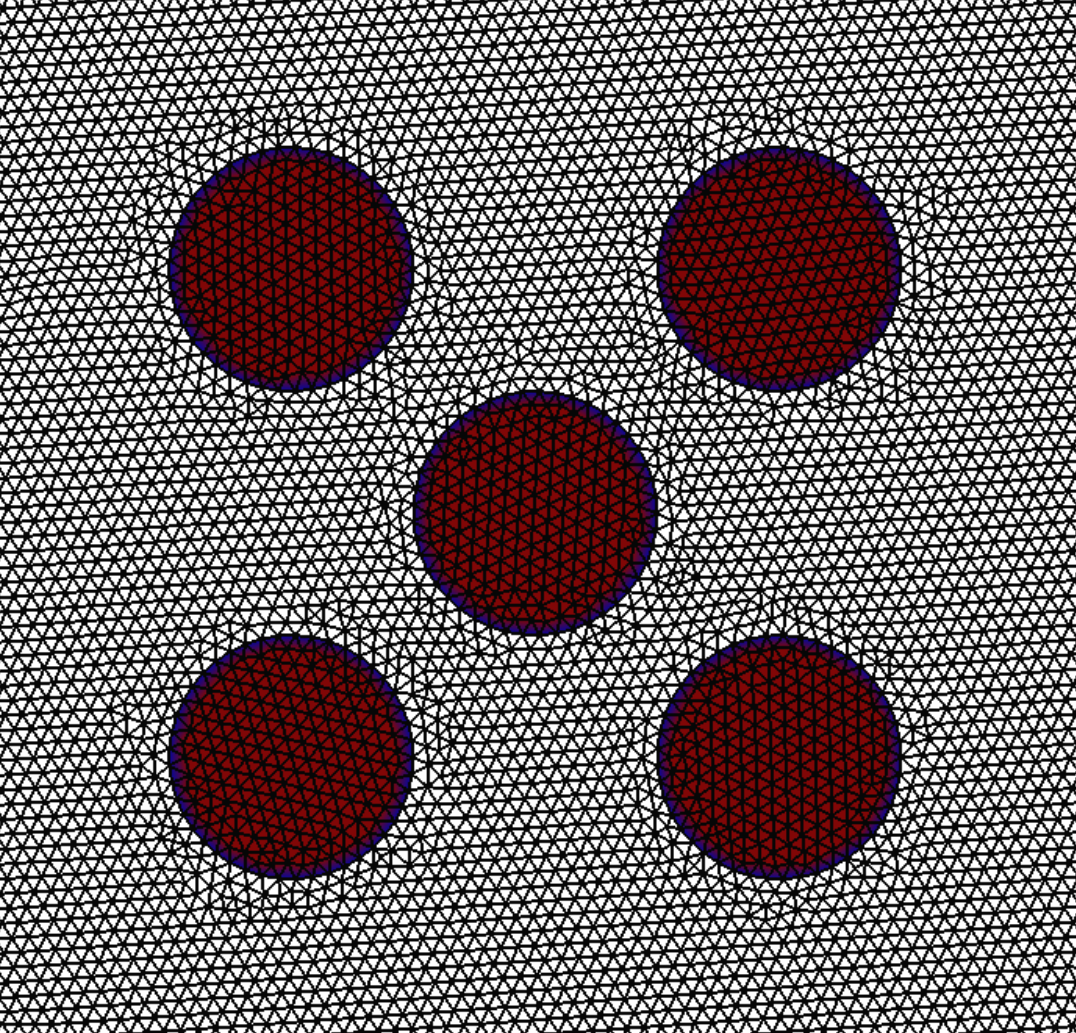}
    \includegraphics[width=0.24\textwidth]{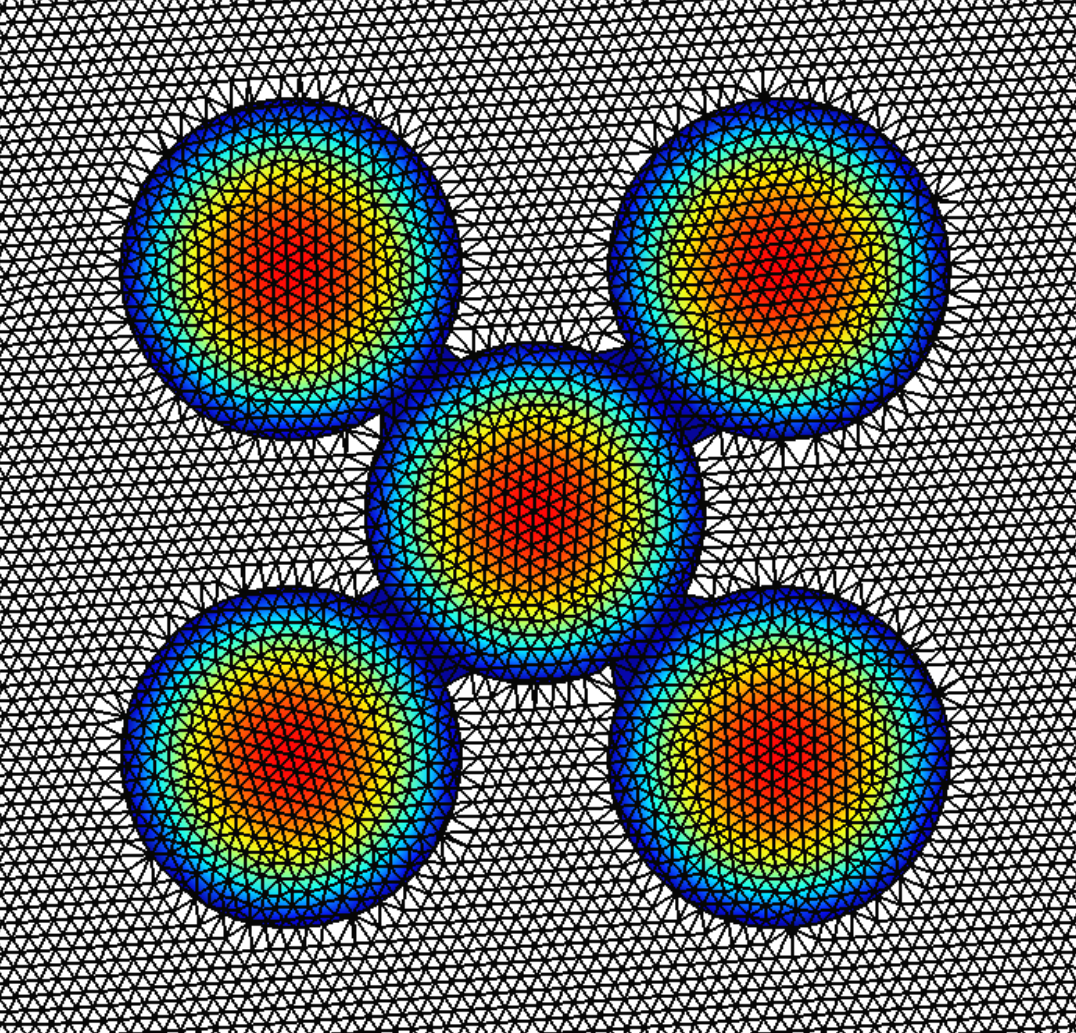}
    \includegraphics[width=0.24\textwidth]{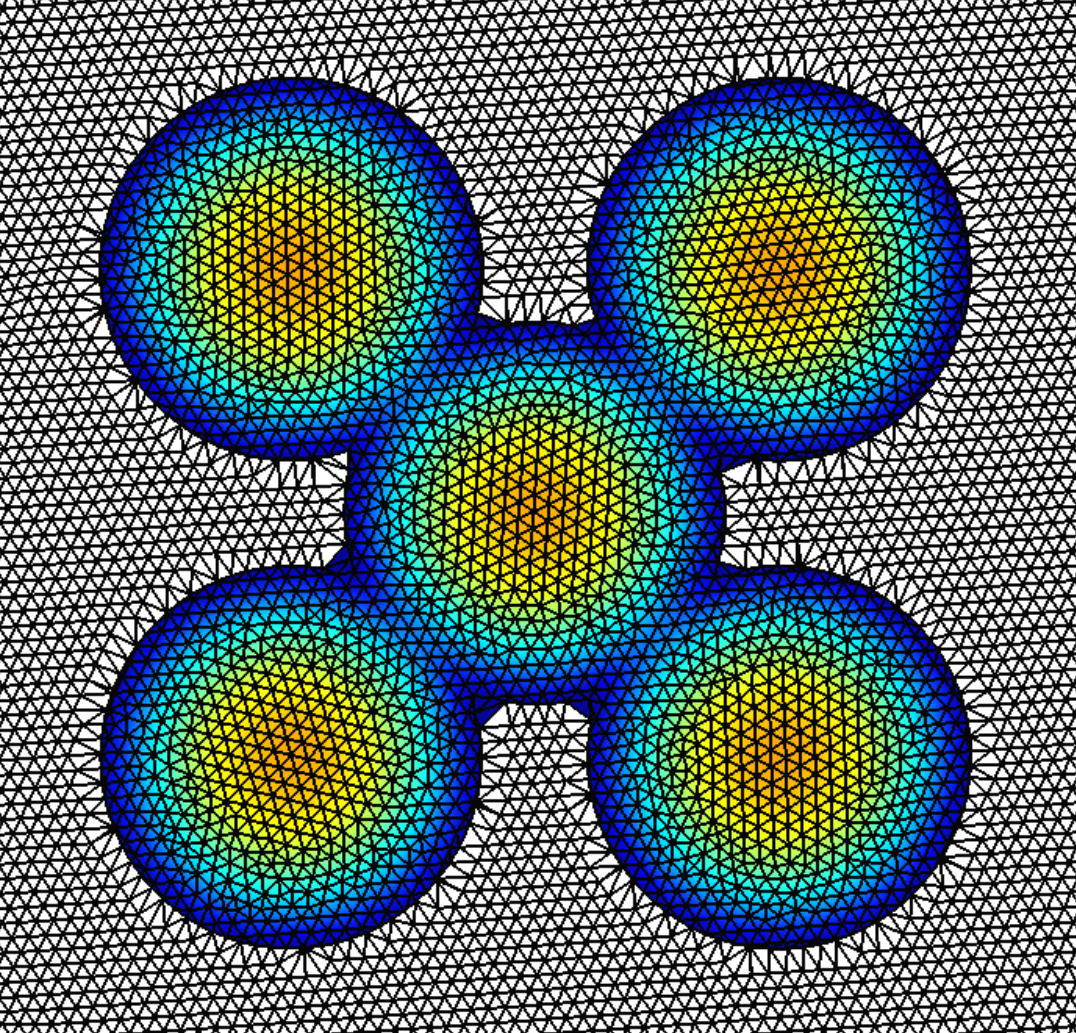}
    \includegraphics[width=0.24\textwidth]{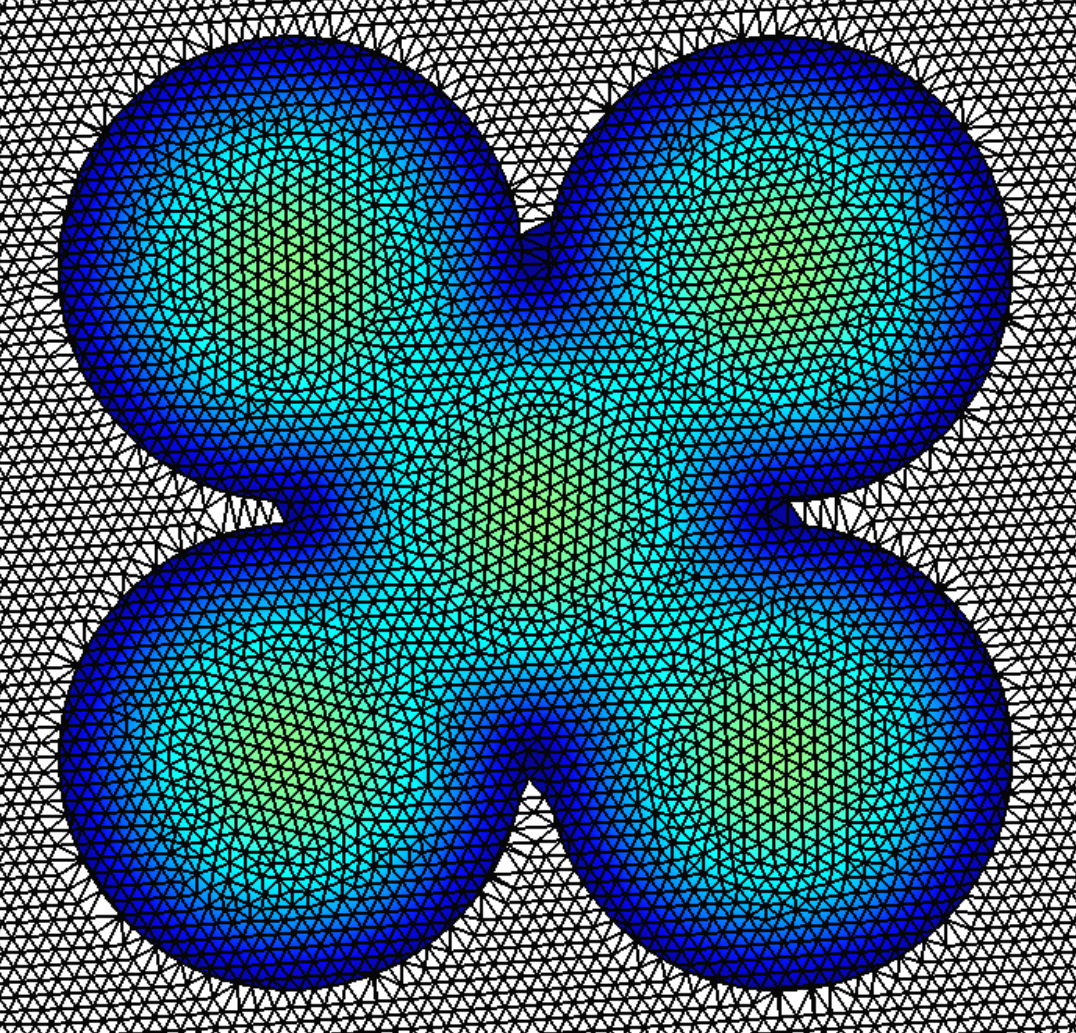}
  \end{center}
  \caption{Four sequential snapshots illustrating the evolution of the solution of the PME associated with the initial condition depicted on the left. Interface coalescence is seamlessly handled by the \xm{} approach.}
  \label{fig:coal}
\end{figure}

The \xm{} method naturally handles cases where the topology of the interface $\Gamma(u(t))$ changes over time through coalescence or the disappearance of parts of the interface. Results obtained with an initial condition $u_0$ whose support is composed of five disjoint disks, are shown in Figure~\ref{fig:coal}. At $t_0$, interface $\Gamma(u_0)$ consists of five non-intersecting circles. Over time, the support of solution $u$ becomes simply connected, resulting in interface $\Gamma$ forming a single closed curve.

\begin{figure}[h!]
  \begin{center}
    \includegraphics[width=0.24\textwidth]{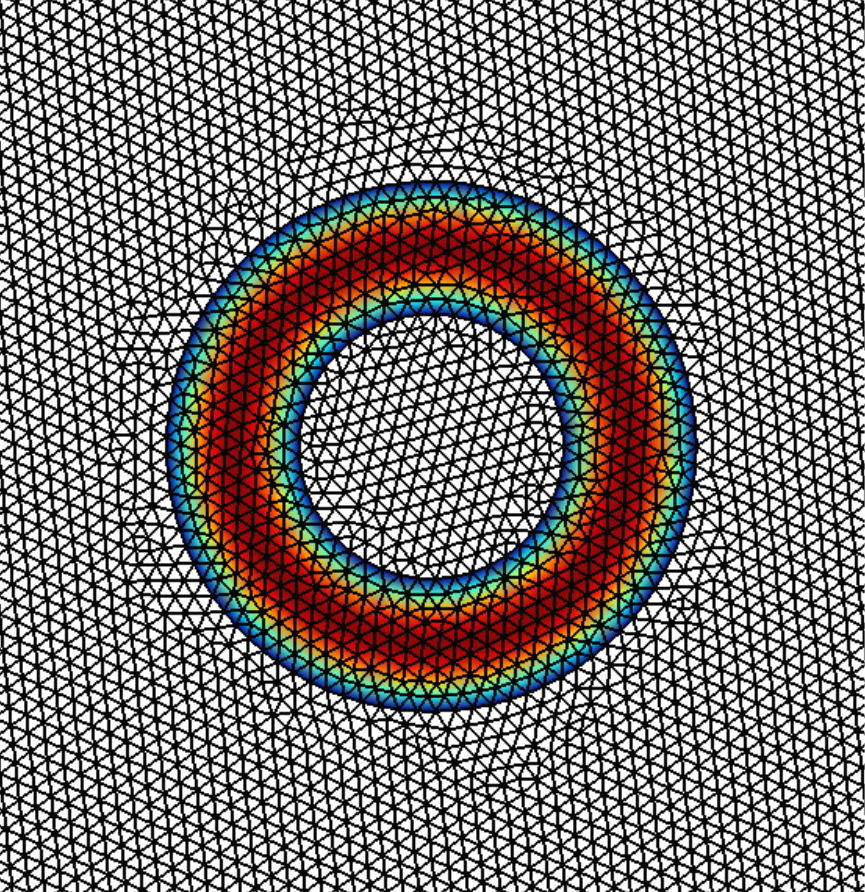}
    \includegraphics[width=0.24\textwidth]{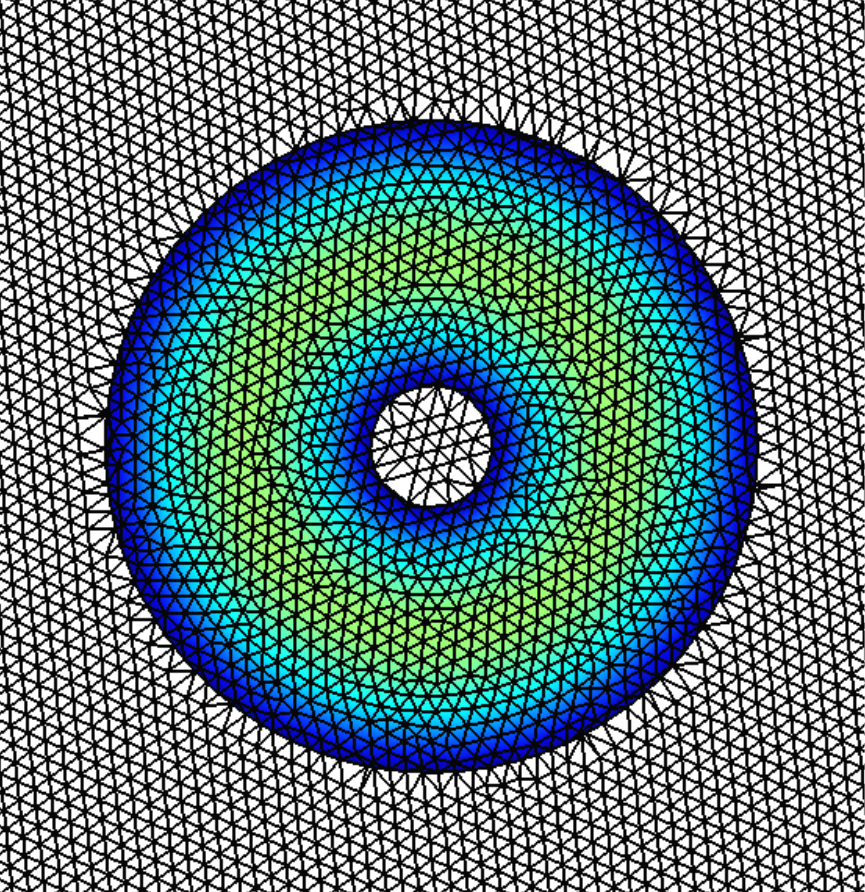}
    \includegraphics[width=0.24\textwidth]{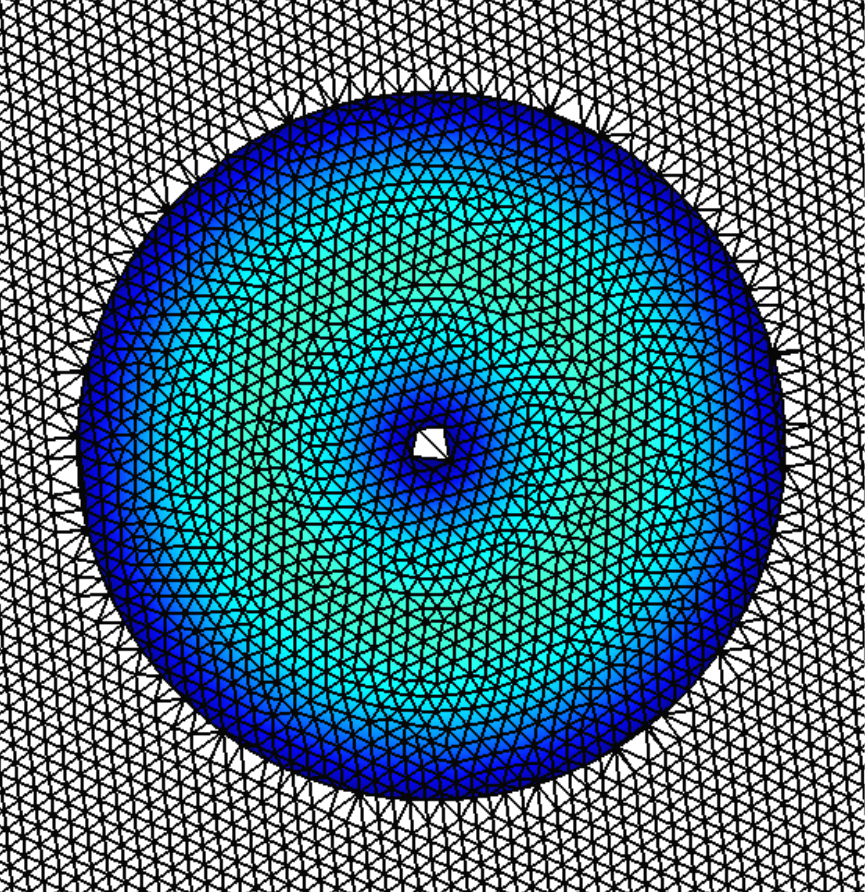}
    \includegraphics[width=0.24\textwidth]{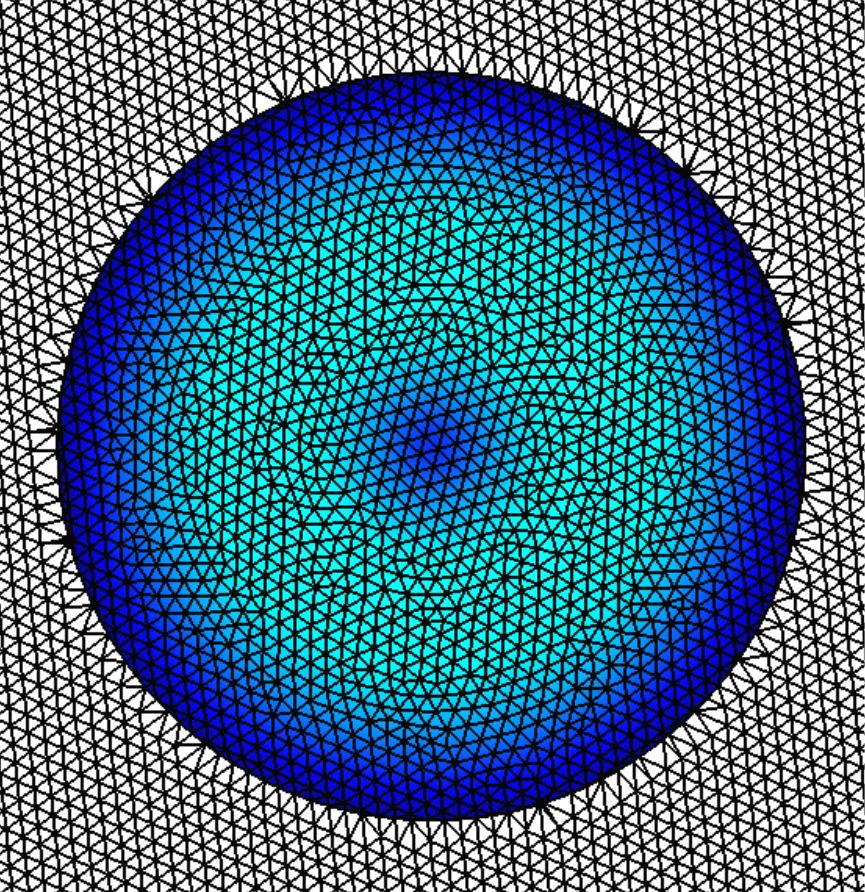}
  \end{center}
  \caption{Four sequential snapshots illustrating the evolution of the solution of the PME where a part of the interface vanishes. The initial condition is depicted on the left.}
  \label{fig:vanish}
\end{figure}

Parts of interface $\Gamma(u(t))$ can also vanish. This can be observed by using an initial condition $u_0$ whose support is an annulus. Initially, interface $\Gamma(u_0)$ comprises two non-intersecting circles. However, at a certain time $t$, the support of $u$ condenses into a disk, causing interface $\Gamma(u(t))$ to become a single circle. Results obtained using the \xm{} approach are depicted in Figure~\ref{fig:vanish}.

These two examples demonstrate that the \xm{} method effectively handles changes in interface topology, accommodating instances of coalescence and vanishing without requiring any additional attention or specific adjustments. The method accurately captures the transition from complex, multi-component interfaces to simpler, single-component configurations.

\section{Conclusion}
In this paper, we have presented the application of the eXtreme Mesh deformation (\xm{}) approach to address the challenges associated with the numerical simulation of the Porous Medium Equation (PME). While the \xm{} approach, designed for tracking sharp interfaces without the need for remeshing or altering mesh topology, has previously found success in solving two-phase problems such as the Stefan phase change model and two-phase flows, our objective here was to extend its applicability to the PME, which constitutes a one-phase problem.

A significant challenge discussed in the literature is the numerical resolution of the non-regularized PME while ensuring non-negativity of solutions, mass conservation, and accurate depiction of interface dynamics. The \xm{} approach, as adapted and implemented in this paper, offers a comprehensive solution to these challenges.

Our numerical experiments served to validate the efficiency of the \xm{} approach in addressing various phenomena inherent to the PME. Through simulations on the Barenblatt-Pattle solution, we demonstrated significant improvements in interface localization precision, modest enhancements in $L^2$ convergence, and the successful preservation of mass conservation. Additionally, our investigation into waiting time phenomena reaffirmed the \xm{} approach's ability to accurately reproduce analytical results, thereby establishing its reliability and accuracy. Furthermore, we showcased the versatility and robustness of the \xm{} approach in handling interface topology changes, demonstrating its seamless adaptability in accommodating instances of coalescence and vanishing without necessitating additional adjustments or specific treatments.

In conclusion, while the \xm{} approach has been previously utilized for Stefan problems and two-phase flows, our application to the PME domain represents a significant extension of its utility. The demonstrated effectiveness of the approach in addressing longstanding challenges and consistently delivering reliable results establishes it as a valuable tool for solving the PME.

\section*{Acknowledgement}
This project has received funding from the European Research Council (ERC) under the European Union’s Horizon Research and Innovation Programme (Grant agreement No. 101 071 255).

\appendix
\newpage
\section{Justification of the mesh velocity definition}
\label{ann:meshVel}

In this appendix is detailed the choice made in Section~\ref{sec:meshVel} for the mesh velocity $\vec{w}$ inside a time interval $[t_n, t_{n+1}]$. As previously mentioned, $\vec{w}$ verifies:
\begin{equation}
  \dint_{t_n}^{t_{n+1}} \vec{w}(\vec{X}_0, t) \,\text{d}t = \vec{X}(\vec{X}_0, t_{n+1}) - \vec{X}(\vec{X}_0, t_{n}) = \vec{X}_{n+1} (\vec{X}_0) - \vec{X}_{n} (\vec{X}_0)\text{, }\forall \vec{X}_0\in\Omega
  \label{eq:constrWann}
\end{equation}
and we have the flexibility to choose any $\vec{w}$ satisfying Eq.~\eqref{eq:constrWann}. The natural modeling choice proposed in \cite{moes2023extreme} is a mesh velocity $\vec{w}$ such as $\left.\dpfrac{\vec w}{t}\right|_{\vec X} = \vec 0$ in $[t_n, t_{n+1}]$. Using definition established in Eq.~\eqref{eq:defXp}, it writes:
\begin{equation}
  \left\{
  \begin{array}{l}
    \vspace{0.5em}
    \vec{w}(\vec X) = \dsum_{i\in\mathcal{N}} \vec W^i \phi_i(\vec{X})\\
    \vec{W}^i = \dfrac{\vec{X}_{n+1}(i) - \vec{X}_{n}(i)}{t_{n+1}-t_n}\\
  \end{array}
  \right.
  \label{eq:defaultWann}
\end{equation}

This option is suitable for addressing \twop\, problems, however it is not appropriate for \onep\, problems due to its inadequate representation of the underlying physics. Indeed, a fully converged solution (there is no Lagrange multipliers left and $R((U_{n}, \vec{X}_{n}), (U_{n+1}, \vec{X}_{n+1})) = 0$) can exhibit an unwanted phenomenon where the phase $\mathcal{P}_u$ spreads more than expected and on a very thin sheet (illustrated Figure~\ref{fig:oldW}). This behavior is the result of the fact that Eq.~\eqref{eq:defaultWann} is unable to model the interface motion inside a time interval.

To illustrate this, we will examine the behavior of this modelization of mesh velocity on a straightforward 1D transport problem:
\begin{equation}
  \left\{
  \begin{array}{lcll}
    \vspace{0.5em}
    \dpfrac{u}{t} & = & -\vec v.\nabla u = -v\dpfrac{u}{x} &\text{, } v\in\mathbb{R}^{+*}\\
    \vspace{0.5em}
    u(x, 0) & = & u_0(x)&\text{, } \forall (x, t)\in\mathbb{R}\times[0, T]\\
    u_0(x) & = &
    \left\{
    \begin{array}{l}
      -ax\\
      -bx\\
    \end{array}
    \right.
    &
    \text{, }
    \begin{array}{l}
      \hspace{-0.4em}\text{if }x<0\text{, } a > 0\\
      \hspace{-0.4em}\text{if }x\geq0\text{, } b\geq 0\\
    \end{array}
    \\
  \end{array}
  \right.
  \label{eq:adv1Dex}
\end{equation}
which admit as solution:
\begin{equation}
  u(x,t) = u_0(x-vt)\text{, }\forall (x,t)\in\mathbb{R}\times[0, T]
\end{equation}

The case $b=0$ corresponds to a \onep\, problem (the \emph{phase} $\mathcal{P}_{u(t)} = \{x\in\mathbb{R}\mid u(x, t)>0\}$ and the interface $\Gamma(u(t)) = \partial\mathcal{P}_{u(t)}$) and the case $b>0$ to a \twop\, problem (the first phase $\mathcal{P}^1_{u(t)}=\{x\in\mathbb{R}\mid u(x,t)>0\}$, the second phase $\mathcal{P}^2_{u(t)}=\{x\in\mathbb{R}\mid u(x,t)<0\}$ and the interface $\Gamma(u(t)) = \partial\mathcal{P}^1_{u(t)} = \partial\mathcal{P}^2_{u(t)} = \{ x\in\mathbb{R} \mid u(x, t) = 0\}$).

We are now applying to Eq.~\eqref{eq:adv1Dex} the space and time discretization detailed Section~\ref{sec:xmesh} for one time interval $[0, \Delta t]$ and with four nodes $\{N_0, N_1, N_2, N_3\}$ and the corresponding three elements $\{[N_0, N_1], [N_1, N_2], [N_2, N_3]\}$. Nodes $\{N_0, N_1, N_3\}$ position will be fixed, node $N_2$ position will be time-dependent and its velocity $W(N_2)$ is constant.
\begin{equation}
  \left\{
  \begin{array}{l}
    \vspace{0.5em}
    X(N_0, t) = -1\\
    X(N_1, t) = 0\\
    X(N_2, t) = 1 + W(N_2)t = 1 + \dfrac{\Delta x}{\Delta t} t\\
    X(N_3, t) = 2\\
  \end{array}\right.
\end{equation}

\begin{figure}[!h]
  \centering
  \includegraphics[width=0.45\textwidth]{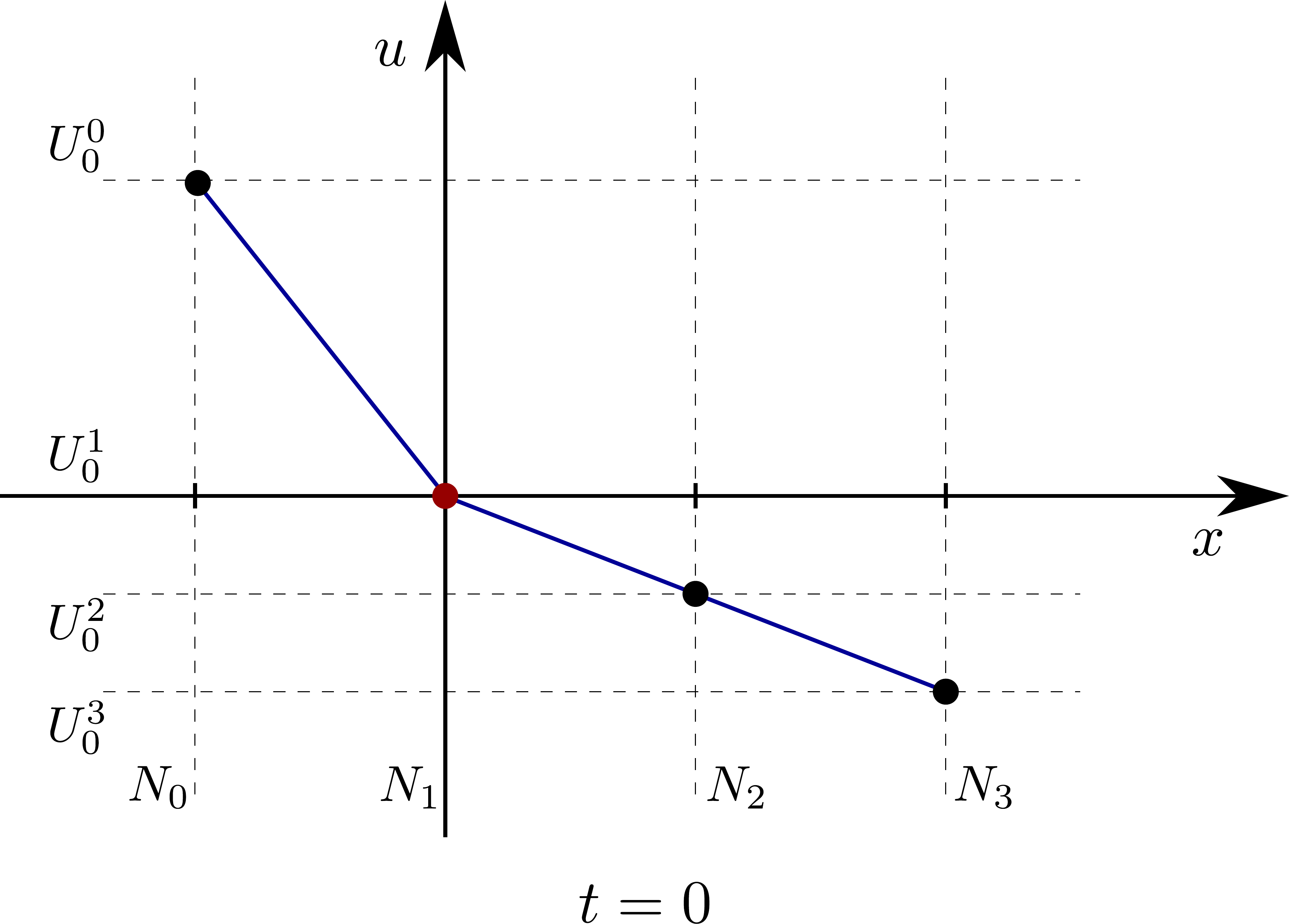}
  \hspace{2em}
  \includegraphics[width=0.45\textwidth]{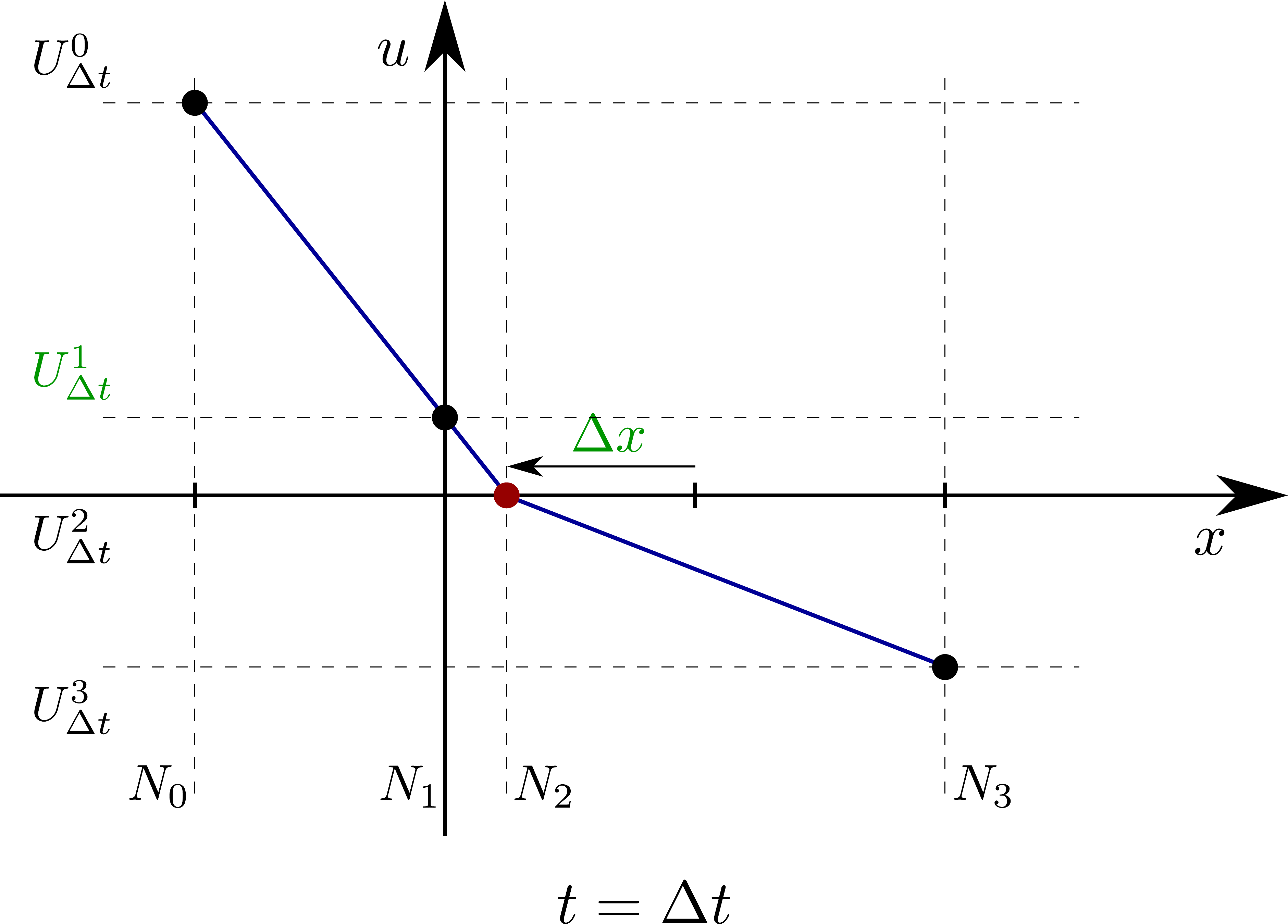}
  \caption{Representation of the discretized solution at $t=0$ and $t=\Delta t$. Red dots are representing interface $\Gamma(u)$ location and quantities marked in green correspond to unknowns of the problem.}
  \label{fig:spaceTimeDiscr}
\end{figure}

Finally, we want node $N_2$ to be located exactly on the interface at $t=\Delta t$, meaning that $U^2_{\Delta t} = 0$. We therefore have the following initial and boundary conditions:
\begin{equation}
  \left\{
  \begin{array}{l}
    U^0_0 = a\\
    U^1_0 = 0\\
    U^2_0 = -b\\
    U^3_0 = -2b\\
  \end{array}\right.,\,
  \left\{
  \begin{array}{l}
    U^0_{\Delta t} = a(1+v\Delta t)\\
    U^2_{\Delta t} = 0\\
    U^3_{\Delta t} = b(-2+v\Delta t)\\
  \end{array}\right.\text{. }
\end{equation}
The two unkowns of the problem are $U^1_{\Delta t}$ and $\Delta x$. The discretization and constraints are summed up Figure~\ref{fig:spaceTimeDiscr}.

\begin{figure}[!h]
  \centering
  \includegraphics[width=0.45\textwidth]{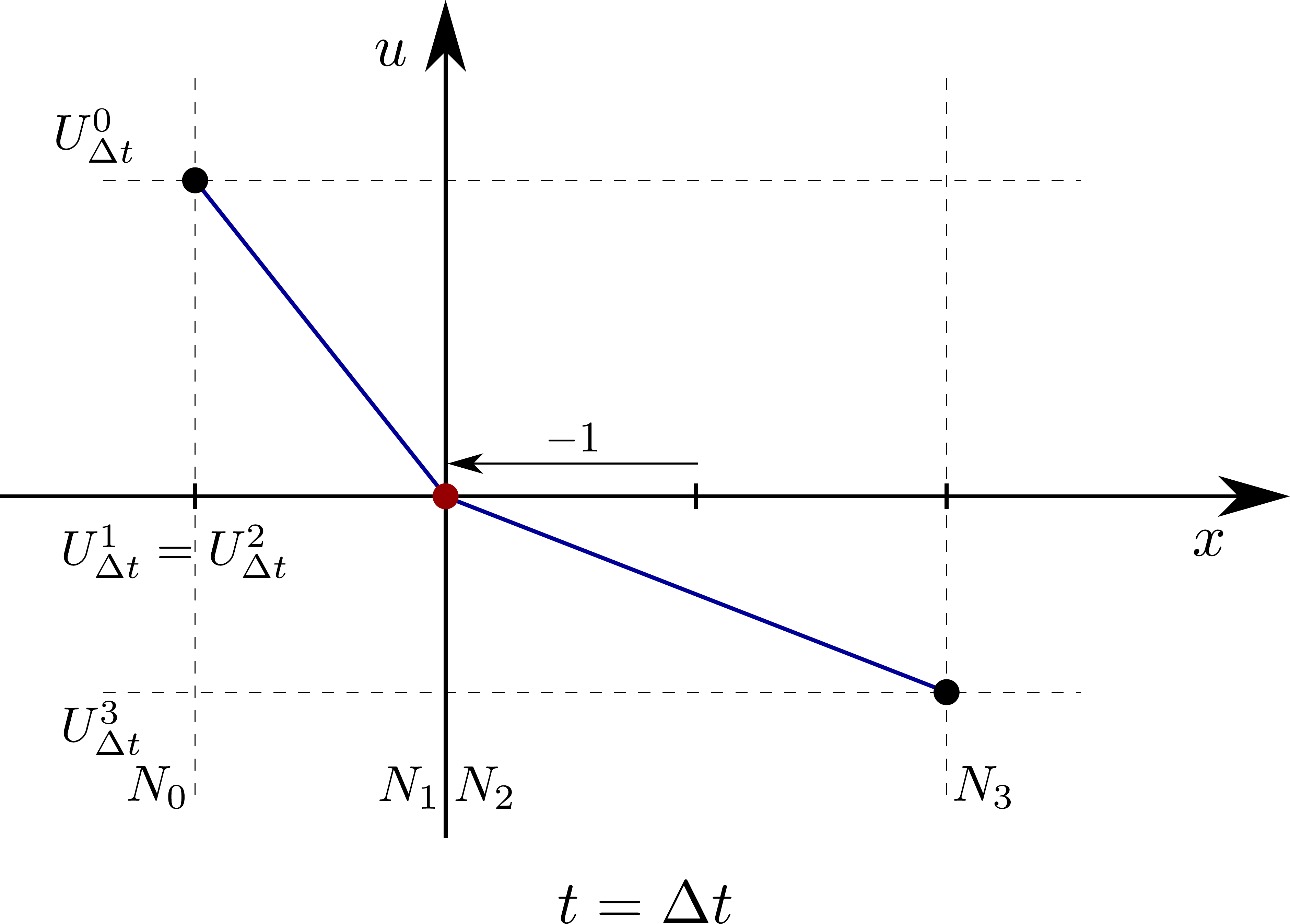}
  \hspace{2em}
  \includegraphics[width=0.45\textwidth]{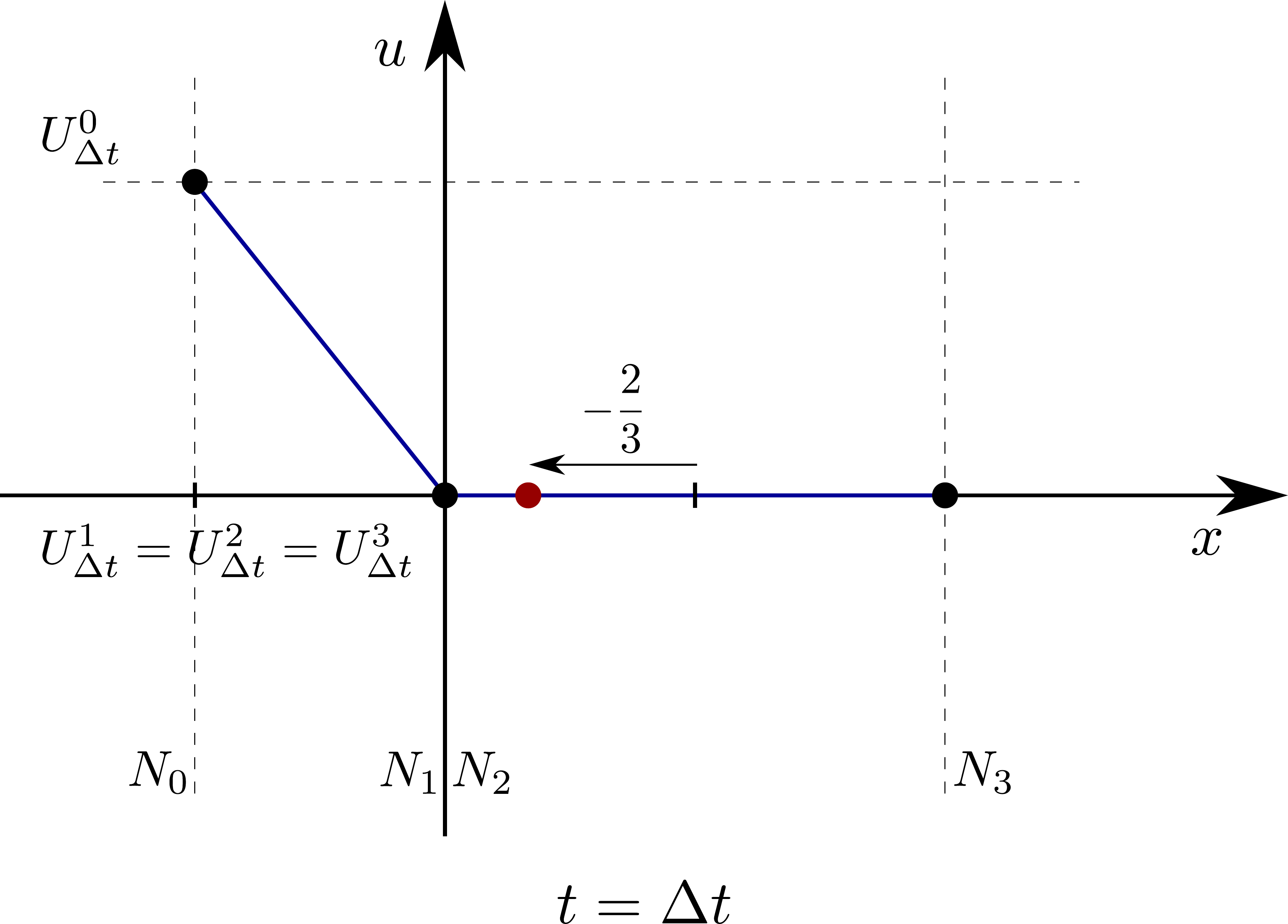}
  \caption{Discretized solutions for the \twop\, problem (left) and the \onep\, problem (right) when $\Delta t \to 0$. The red dot represent the interface location and the arrow the node $N_2$ displacement.}
  \label{fig:advIntdt0}
\end{figure}

The discretized problem admits the following solution:
\begin{equation}
  \begin{array}{ll}
    \vspace{0.5em}
    \text{If } b=0,\,&
    \left\{
    \begin{array}{lcl}
      \vspace{0.5em}
      U^1_{\Delta t} & = & \dfrac{av\Delta t(3v\Delta t+4)}{3v\Delta t + 6}\\
      \Delta x & = & \dfrac{3v\Delta t - 2}{3}\\
    \end{array}
    \right.\\
    \text{If } b>0,\,&
    \left\{
    \begin{array}{lcl}
      \vspace{0.5em}
      U^1_{\Delta t} & = & \dfrac{16a+5b}{21}v\Delta t + o(\Delta t)\\
      \Delta x & = & -1 + \dfrac{2a + 19b}{21b}v\Delta t + o(\Delta t)\\
    \end{array}
    \right.\\
  \end{array}
  \label{eq:solAdv1Ddiscr}
\end{equation}
We can notice that for the \twop\, problem (case $b>0$),  $\displaystyle\lim_{\Delta t \to 0} \Delta x = -1$  ($N_2$ tends towards $N_1$ when $\Delta t\rightarrow 0$) which is the expected behavior. However, for the \onep\, problem (case $b=0$) we have $\displaystyle\lim_{\Delta t \to 0} \Delta x = -\frac{2}{3}$, exhibiting the pathological behavior where the phase $\mathcal{P}_u$ spreads more than expected on a thin sheet. These solutions are depicted Figure~\ref{fig:advIntdt0}.

\begin{figure}[!h]
  \centering
  \includegraphics[width=1.0\textwidth]{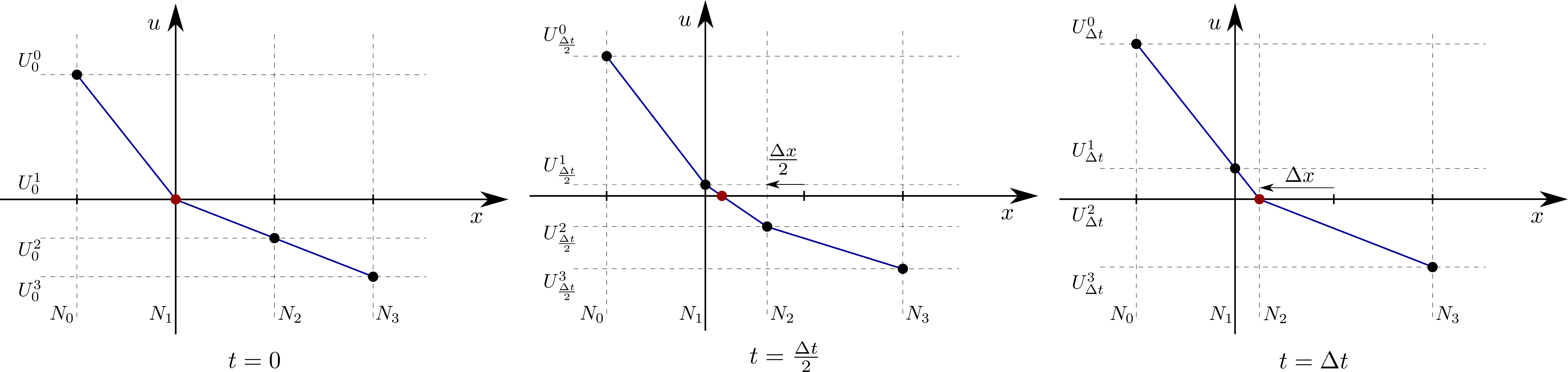}
  \caption{Behavior of the discretized solution of the \twop\, problem in the time interval $[0, \Delta t]$. The red dot marks the location of the interface ($u^h(x, t)=0$). The interface moves from left to right as expected.}
  \label{fig:behavSol2P}
\end{figure}

\begin{figure}[!h]
  \centering
  \includegraphics[width=1.0\textwidth]{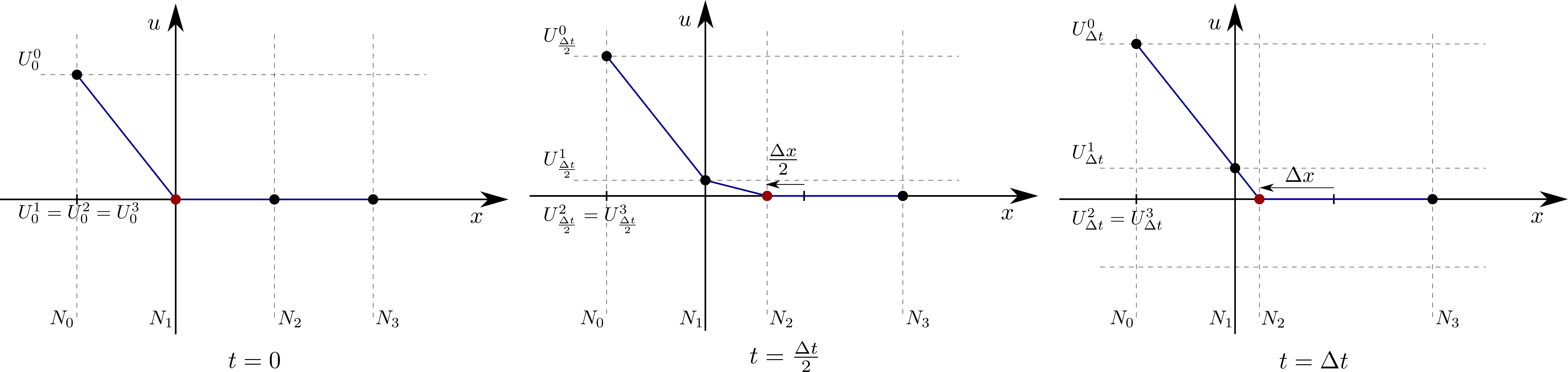}
  \caption{Behavior of the discretized solution for the \onep\, problem in the time interval $[0, \Delta t]$. The red dot marks the location of the interface ($u^h(x, t)=0$). The interface moves from right to left. This behavior, imposed by the discretization choices made, is the opposite of what we want to model.}
  \label{fig:behavSol1P}
\end{figure}

The space discretization performed, the time integration scheme chosen and the mesh velocity modelization choice made are implying in particular that $N_2$ varies continuously from $N_2(0)$ to $N_2(\Delta t)$ inside the time interval, and that $U^i$ varies continuously from $U^i(0)$ to $U^i(\Delta t)$. A depiction of the solution behavior inside the time interval $[0, \Delta t ]$ can be found Figure~\ref{fig:behavSol2P} for the \twop\, problem and Figure~\ref{fig:behavSol1P} for the \onep\, problem.

The important behavior to notice is that for the \twop\, problem, the discretized solution inside the time interval exhibits an interface moving from the left to the right as expected, but in the case of the \onep\, problem the discretized solution inside the time interval exhibits an interface moving from the right to the left. This behavior appears if the node chose to be moved on the interface $\Gamma^h_{n+1}$ was belonging to the empty region (defined section \ref{sec:PME}) at $t = t_n$. In other words this behavior occurs if $p\in\mathcal{Q}_{u^h_n}$. The mesh velocity model described by Eq.~\eqref{eq:defaultWann} fails to represent the behavior of the interface for the \onep\, problem and therefore requires modification.

\begin{figure}[!h]
  \centering
  \includegraphics[width=1.0\textwidth]{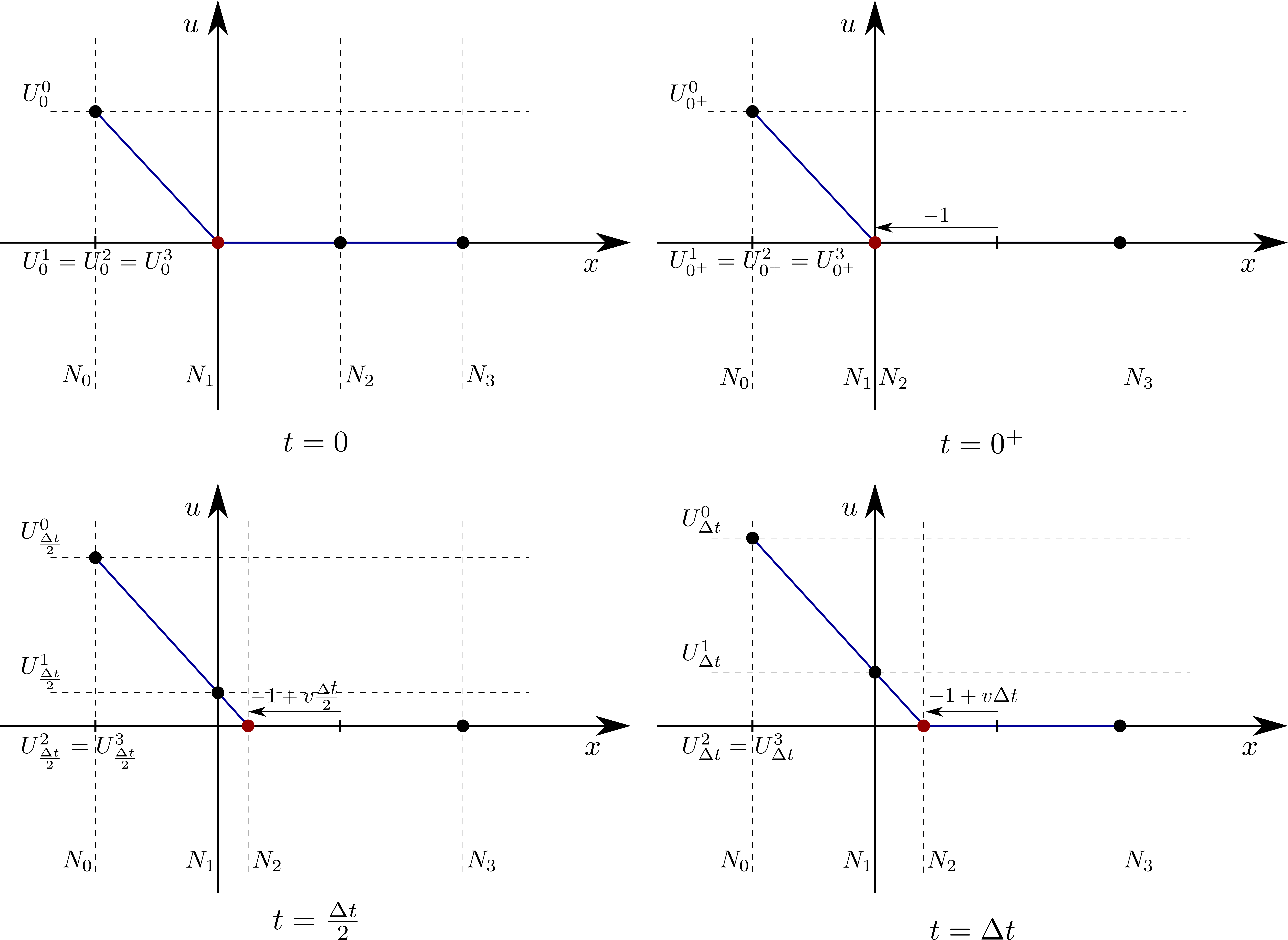}
  \caption{Behavior of the discretized solution for the problem in the time interval $[0, \Delta t]$ when the mesh velocity is defined with Eq.~\eqref{eq:newWdefSimple}. The red dot marks the location of the interface ($u^h(x, t)=0$). The interface now moves from left to right, which is consistent with the physics we want to model.}
  \label{fig:behavSol1Pcorr}
\end{figure}

The choice we are making for $\vec w$ is:
\begin{equation}
  \left\{
  \begin{array}{ll}
    \vec{W}^p(t) = (\vec{X}_{n}(q_p^0) - \vec{X}_{n}(p))\delta_{t_n}(t) + \dfrac{\vec{X}_{n+1}(p) - \vec{X}_{n}(q_p^0)}{t_{n+1}-t_n}&\text{, if $p\in\mathcal{Q}_{u^h_n}$}\\
    & \text{and $q_p^0$ defined Eq.~\eqref{eq:qp}}\\
    \vec{W}^p(t) = \dfrac{\vec{X}_{n+1}(p) - \vec{X}_{n}(p)}{t_{n+1}-t_n}&\text{, else.}\\
  \end{array}
  \right.
  \label{eq:defWbase}
\end{equation}
where $\delta_{t_n}$ is a Dirac distribution in time centered in $t_n$. This expression for $\vec w$ can be simplified when taking into account that the problem is a \onep\, problem. Indeed, the weak formulation of Eq.~\eqref{eq:adv1Dex} discretized in space writes:
\begin{equation}
  \begin{array}{lll}
    \dint_{t_n}^{t_{n+1}}\dfrac{\text{d}}{\text{d}t} \dint_\Omega \bar{u}^h u^h \id\Omega\id t & = & - \dint_{t_n}^{t_{n+1}}\dint_\Omega \bar{u}^h \vec v.\nabla u^h \id\Omega\id t - \dint_{t_n}^{t_{n+1}}\dint_\Omega u^h\vec{w}\nabla \bar{u}^h\id\Omega\id t.\\
  \end{array}
  \label{eq:weakAdv1Dex}
\end{equation}
The term of Eq.~\eqref{eq:weakAdv1Dex} involving the mesh speed $\vec w$ can be developed as:
\begin{equation}
  \begin{array}{lll}
    \dint_{t_n}^{t_{n+1}}\dint_\Omega u^h\vec{w}\nabla \bar{u}^h\id\Omega\id t & = & \dint_{t_n}^{t_{n+1}}\dint_\Omega u^h\dsum_{i\in\mathcal{N}} \phi_i(\vec{X})\vec W^i(t).\nabla \bar{u}^h\id\Omega\id t\\
    & = & \dsum_{i\in\mathcal{N}}\dint_\Omega \dint_{t_n}^{t_{n+1}} u^h \phi_i(\vec{X})\vec W^i(t).\nabla \bar{u}^h\id t\id\Omega\\
  \end{array}
  \label{eq:proofW1}
\end{equation}
Examining the term of the sum corresponding to the node $i=p\in\mathcal{Q}_{u^h_n}$ and utilizing expression Eq.~\eqref{eq:defWbase}, we derive:
\begin{equation}
  \begin{array}{lll}
    \dint_\Omega \dint_{t_n}^{t_{n+1}} u^h \phi_p(\vec{X})\vec W^p(t).\nabla \bar{u}^h\id t\id\Omega & = & \dint_\Omega u^h_n \phi_p(\vec{X})(\vec{X}_{n}(q_p^0) - \vec{X}_{n}(p)).\nabla \bar{u}^h\id\Omega \\
    & + & \dint_\Omega \dint_{t_n}^{t_{n+1}} u^h \phi_p(\vec{X})\dfrac{\vec{X}_{n+1}(p) - \vec{X}_{n}(q_p^0)}{t_{n+1}-t_n}.\nabla \bar{u}^h\id t\id\Omega\\
  \end{array}
  \label{eq:proofW2}
\end{equation}
As $p\in\mathcal{Q}_{u^h_n}$, we know that:
\begin{equation}
  u^h(\vec{X})\phi_p(\vec{X}) = 0\text{, } \forall\vec X\in\Omega
  \label{eq:nullUphi}
\end{equation}
Using Eq.~\eqref{eq:nullUphi} in Eq.~\eqref{eq:proofW2} it remains:
\begin{equation}
    \dint_\Omega \dint_{t_n}^{t_{n+1}} u^h \phi_p(\vec{X})\vec W^p(t).\nabla \bar{u}^h\id t\id\Omega = \dint_\Omega \dint_{t_n}^{t_{n+1}} u^h \phi_p(\vec{X})\dfrac{\vec{X}_{n+1}(p) - \vec{X}_{n}(q_p^0)}{t_{n+1}-t_n}.\nabla \bar{u}^h\id t\id\Omega
  \label{eq:proofW3}
\end{equation}
The choice for $\vec w$ described Eq.~\eqref{eq:defWbase} is therefore equivalent to the simplier:
\begin{equation}
  \left\{
  \begin{array}{l}
    \vec{W}^p(t) = \dfrac{\vec{X}_{n+1}(p) - \vec{X}_{n}(q_p^0)}{t_{n+1}-t_n}\text{, if $p\in\mathcal{Q}_{u^h_n}$}\\
    \vec{W}^p(t) = \dfrac{\vec{X}_{n+1}(p) - \vec{X}_{n}(p)}{t_{n+1}-t_n}\text{, else.}\\
  \end{array}
  \right.
  \label{eq:newWdefSimple}
\end{equation}
Using this expression for $\vec w$, we now obtain the behavior of the solution inside the time interval illustrated Figure~\ref{fig:behavSol1Pcorr}.

This behavior is now properly modeling evolution of interface $\Gamma$ position regarding time. When applying this new modelization of $\vec w$, the discretization of problem \ref{eq:adv1Dex} admits for solution:
\begin{equation}
  \text{If } b=0,\,
  \left\{
  \begin{array}{lcl}
    U^1_{\Delta t} & = & av\Delta t\\
    \Delta x & = & -1 + v\Delta t\\
  \end{array}
  \right.
  \label{eq:resWuncorr2}
\end{equation}
for which $\displaystyle\lim_{\Delta t \to 0} \Delta x = -1$, corresponding to $N_2$ tending towards $N_1$ when $\Delta t\rightarrow 0$. $\mathcal{P}_u$ does not spread more than expected anymore.

\newpage
\bibliographystyle{plain}
\bibliography{biblio}

\end{document}